\documentclass[11pt,a4paper]{article}

\usepackage{url,amsmath,amssymb,latexsym,pstricks,mathrsfs,comment,amsthm,graphicx,tikz,tikz-cd,accents,pgffor,cite,wrapfig,multicol,float,cases,calc,bibspacing,geometry,bbm,arydshln,enumerate}
\usepackage[colorlinks]{hyperref}
\usepackage[T1]{fontenc}
\usepackage{ wasysym }
\usepackage{stmaryrd}
\usepackage{enumitem}
\usetikzlibrary{calc}

\numberwithin{equation}{section}

\newtheorem{thm}[equation]{Theorem}
\newtheorem{lemma}[equation]{Lemma}

\newtheorem{prop}[equation]{Proposition}

\theoremstyle{definition}

\newtheorem{rem}[equation]{Remark}

\newcommand{\nc}{\newcommand}
\nc{\rnc}{\renewcommand}

\let\oldproofname=\proofname
\rnc{\proofname}{\rm\bf{\oldproofname}}

\rnc{\arraystretch}{1.2}

\rnc\L{\mathscr L}
\nc\hs\heartsuit
\nc\ds\diamondsuit

\nc\TT{\mathcal T}
\rnc\S{\mathcal S}
\nc\Inj{\operatorname{Inj}}
\nc\Surj{\operatorname{Surj}}

\rnc\ss\spadesuit
\nc\Fup{\Fu^+}
\nc\Lp{\L^+}

\nc\Tone{\textup{\textsf{T$_{\text{1}}$}}}
\nc\Ttwo{\textup{\textsf{T$_{\text{2}}$}}}
\nc\Tthree{\textup{\textsf{T$_{\text{3}}$}}}
\nc\Ttwod{\textup{\textsf{T$_{\text{2}}'$}}}
\nc\Tthreed{\textup{\textsf{T$_{\text{3}}'$}}}
\nc\Tfour{\textup{\textsf{T$_{\text{4}}$}}}
\nc\Ti{\textup{\textsf{T$_{\text{i}}$}}}
\nc{\T}{\textup{\textsf{T}}}

\nc\ben{\begin{enumerate}[label=\textup{(\roman*)},leftmargin=7mm]}
\nc\benT{\begin{enumerate}[label=\textsf{(T$_{\text{\arabic*}}$)},leftmargin=10mm]}
\nc\benTd{\begin{enumerate}[label=\textsf{(T$_{\text{\arabic*}}'$)},leftmargin=10mm]}
\nc\een{\end{enumerate}}

\nc\Mon{\mathcal M}
\nc\M{\mathcal M}
\nc\Fu{\mathscr F}
\nc\leqJ{\leq_{\mathscr J}}
\nc\geqJ{\geq_{\mathscr J}}
\nc\Sub{\operatorname{Sub}}
\nc\Pos{{\bf P}}

\nc\gL{\mathrel{\mathscr L}}
\nc\gR{\mathrel{\mathscr R}}
\nc\gH{\mathrel{\mathscr H}}
\nc\gJ{\mathrel{\mathscr J}}
\nc\gD{\mathrel{\mathscr D}}
\nc\gK{\mathrel{\mathscr K}}

\nc{\Op}{\mathcal O}

\rnc{\O}{\mathbb O}
\nc{\G}{\mathbb G}
\nc{\E}{\mathbb E}
\nc{\F}{\mathbb F}
\rnc{\P}{\mathbb P} 
\nc{\I}{\mathbb I}
\nc{\Q}{\mathbb Q}
\nc\GL{\G_L}
\nc\GR{\G_R}
\nc\GLR{\G_{LR}}
\nc\FL{\F_L}
\nc\FR{\F_R}
\nc\FLR{\F_{LR}}
\nc\PLR{\P_{LR}}
\nc\PL{\P_L}
\nc\PR{\P_R}
\nc{\U}{\mathbb U}
\nc{\V}{\mathbb V}
\nc{\W}{\mathbb W}

\nc\GH{\G_\heartsuit}
\nc\PH{\P_\heartsuit}
\nc\FH{\F_\heartsuit}

\nc{\PB}{\mathcal{PB}} 
\nc{\X}{\mathbb X}
\nc{\Y}{\mathbb Y}
\nc{\Z}{\mathbb Z}

\nc{\set}[2]{\{ {#1} : {#2} \}} 

\rnc{\emptyset}{\varnothing}
\nc{\sm}{\setminus}
\nc{\sub}{\subseteq}
\nc{\la}{\langle}
\nc{\ra}{\rangle}

\nc{\OR}{\qquad\text{or}\qquad}
\nc{\COMMA}{,\qquad}
\nc{\COMMa}{,\quad}
\nc{\AND}{\qquad\text{and}\qquad}
\nc{\ANd}{\quad\text{and}\quad}
\nc{\ANDSIM}{\qquad\text{and similarly}\qquad}
\rnc{\iff}{\ \Leftrightarrow\ }
\rnc{\implies}{\ \Rightarrow\ }

\nc{\bit}{\begin{itemize}}
\nc{\eit}{\end{itemize}}
\nc{\bmc}{\begin{multicols}}
\nc{\emc}{\end{multicols}}
\nc{\firstpfitem}[1]{#1.}
\nc{\pfitem}[1]{\medskip \noindent #1.}
\nc{\afterpfitem}{\medskip \noindent }
\nc{\pfcase}[1]{\medskip \noindent {\bf Case #1.}}
\nc{\pf}{\begin{proof}}
\nc{\epf}{\end{proof}}
\nc{\epfres}{\hfill\qed}
\nc{\epfreseq}{\tag*{\qed}}

\begin{document}

\title{Idempotents and one-sided units: \\  Lattice invariants and a semigroup of functors on the category of monoids}
\author{James East\\
{\it\small Centre for Research in Mathematics; School of Computing, Engineering and Mathematics,}\\
{\it\small Western Sydney University, Locked Bag 1797, Penrith NSW 2751, Australia.}\\
{\tt\small J.East\,@\,WesternSydney.edu.au}}
\date{}

\maketitle

\begin{abstract}
For a monoid $M$, we denote by $\mathbb G(M)$ the group of units, $\mathbb E(M)$ the submonoid generated by the idempotents, and $\mathbb G_L(M)$ and $\mathbb G_R(M)$ the submonoids consisting of all left or right units.  Writing $\mathcal M$ for the (monoidal) category of monoids, $\mathbb G$, $\mathbb E$, $\mathbb G_L$ and $\mathbb G_R$ are all (monoidal) functors $\mathcal M\to\mathcal M$.  There are other natural functors associated to submonoids generated by combinations of idempotents and one- or two-sided units.  The above functors generate a monoid with composition as its operation.  We show that this monoid has size $15$, and describe its algebraic structure.  We also show how to associate certain lattice invariants to a monoid, and classify the lattices that arise in this fashion.  A number of examples are discussed throughout, some of which are essential for the proofs of the main theoretical results.


\end{abstract}

\tableofcontents

\section{Introduction}\label{s:intro}

Idempotent-generated semigroups arise naturally in many settings, and include semigroups of singular transformations, matrices, partitions, and endomorphisms of various structures \cite{Gray2007,MM2007,Howie1966,JE_IBM,East2011,EG2017,DEG2017, DE5,DE4,DE2016,DEM2016,EF2012,Howie1978,HM1990,BFF1997,BF1995, FG2007,FL1992,Fountain1991,FL1993,Erdos1967,HHR1998, Putcha2006,Putcha1988}.  Free idempotent-generated semigroups associated to abstract biordered sets have long been a crucial tool in the structure theory of (regular) semigroups \cite{Nambooripad1979,Easdown1985,DGR2017,DR2013,GR2012_2,GR2012}.  Many well-known monoids are generated by their idempotents and units, and several studies have calculated the submonoids generated by all idempotents and units of important monoids \cite{EF2012,FL1993,HHR1998,JE_IBM,FL1998,HRH1998,   East2012,CH1974,FitzGerald2010,KM2006,JEgrpm,McAlister1998,FitzGerald2003,EF2013,EF2010}.  

One-sided units have also played an important role in many of the above studies (and others), sometimes implicitly.  As an example, consider the full transformation monoid over a set $X$; this monoid is denoted~$\TT_X$, and consists of all mappings $X\to X$ under composition.  
If functions are composed right-to-left, then the left and right units of $\TT_X$ are precisely the injective and surjective mappings, respectively, while the two-sided units are the bijections, which together form the symmetric group $\S_X$.  It was shown in \cite{HRH1998} that finite~$\TT_X$ is generated by~$\S_X$ and a single idempotent, and that infinite $\TT_X$ is generated by $\S_X$ together with two additional mappings, one a left unit and the other a right unit (both of a certain special form).  It follows that infinite $\TT_X$ is generated by its one-sided units, a property that holds for a number of other important monoids \cite{East2014,AM07,HHMR2003,JE_IBM,East2017,HHR1998}.  The submonoids of $\TT_X$ consisting of all left units or all right units (i.e., all injective or all surjective mappings $X\to X$) have of course been studied in a number of settings as well \cite{East2017,DG2013,MP2011,Mesyan2012,HHR1998,Sutov1961_2}.

To the author's knowledge, the article \cite{JE_IBM} was the first to systematically study submonoids generated by all combinations of idempotents and one- or two-sided units of a monoid, though the article \cite{HHR1998} is a forerunner, as it considered set products of subsets consisting of such elements (in the monoid of partial transformations of an infinite set).  The principal object of study in \cite{JE_IBM} was the so-called partial Brauer monoid $\PB_X$ \cite{Maz1998}, which consists of certain graphs under a natural diagramatic multiplication.  The main results in \cite{JE_IBM} were descriptions of the various submonoids of~$\PB_X$, as well as the relationships between them.  These relationships were described using notions such as relative rank \cite{HRH1998,HHR1998}, Sierpi\'nski rank \cite{Sierpinski1935,Banach1935,MP2012} and the Bergman property \cite{MMR2009,Bergman2006}.  The article \cite{JE_IBM} also contained the beginnings of a general theory of submonoids generated by idempotents and one- or two-sided units of arbitrary monoids.  The current article develops this general theory further, as we now describe.

In Section \ref{s:monoids} we define the submonoids we will be concerned with, and show that the set of all such submonoids of a given monoid~$M$ forms a lattice $\L(M)$; this lattice is a natural invariant of the isomorphism class of~$M$, and its generic shape is shown in Figure~\ref{f:lattice}.  We also show that each such submonoid arises from a functor on the category of monoids, and end the section with a number of examples.  Section \ref{s:prelim} contains preliminary results, mostly concerning intersections of the submonoids, and collapse within the lattice~$\L(M)$; we will also describe some connections with Green's relations and stability of the identity element.  In Section~\ref{s:class} we classify the lattice invariants $\L(M)$, and show that their structure is completely determined by a certain binary quadruple $\T(M)$, which we call the type of $M$.  The main results of Section \ref{s:class} are summarised in Theorem~\ref{t:L}, and the possible shapes of $\L(M)$ are shown in Figures \ref{f:L1}, \ref{f:L2} and \ref{f:L3}.  In Section \ref{s:F} we study the monoid $\Fup$ generated by all of the above-mentioned functors.  This involves calculating all compositions of the functors, and introducing four new ones; in the end we are able to calculate the size of~$\Fup$ and describe its algebraic structure (see Table \ref{t:op2} and Figure \ref{f:div}), using GAP \cite{GAP} for some computations.  Finally, in Section~\ref{s:enhanced} we show that the monoid $\Fup$ may be used to associate a (sometimes) larger lattice $\Lp(M)$ to an arbitrary monoid $M$; we classify these lattices as well (see Figures \ref{f:L4} and \ref{f:L5}), and show that they provide essentially the same information as the original invariant $\L(M)$.

The author would like to acknowledge some valuable comments and questions from a number of colleagues, particularly Robert McDougal, Nik Ru\v skuc, Finn Smith, Timothy Stokes and Lauren Thornton.  The idea to consider the semigroup $\Fup$ traces back to conversations with Dr Thornton about her work on semigroups of operators on radical classes of rings and algebras \cite{Thornton2018,MT2018}.  The author also thanks the referee for their helpful suggestions.

\section{Definitions and basic examples}\label{s:monoids}

In this section we introduce the submonoids (Section \ref{ss:sub}), functors (Section \ref{ss:fun}) and lattices (Section \ref{ss:lat}) that will be at the heart of our investigations, and consider some fundamental examples (Section \ref{ss:egs}).

\subsection{Submonoids}\label{ss:sub}

A monoid is a set $M$ with an associative binary operation, and an identity element $1_M$; the latter will usually be abbreviated to $1$, and the product represented as juxtaposition.  Note that the identity is part of the signature of a monoid, so submonoids must contain the identity, and monoid homomorphisms must map the identity to the identity.

Following the terminology of \cite[Section 1.7]{CPbook}, an element $x$ of a monoid $M$ is:
\bit
\item an \emph{idempotent} if $x=x^2$,
\item a \emph{left unit} if $ax=1$ for some $a\in M$; the element $a$ is a \emph{left inverse} of $x$,
\item a \emph{right unit} if $xa=1$ for some $a\in M$; the element $a$ is a \emph{right inverse} of $x$,
\item a \emph{(two-sided) unit} if it is a left and right unit.
\eit
In general, $x$ could have multiple left or right inverses; however, if it has at least one of each, then it has a unique left unit and a unique right unit, which must be equal, and which we denote by~$x^{-1}$.  
We write
\[
E(M) \COMMA \GL(M) \COMMA \GR(M) \COMMA \G(M)=\GL(M)\cap\GR(M)
\]
for the sets of all idempotents, left units, right units and (two-sided) units of $M$, respectively.  Note that $\GL(M)$, $\GR(M)$ and $\G(M)$ are all submonoids of $M$, with $\G(M)$ a group. We also denote by
\[
\E(M)=\la E(M)\ra
\]
the submonoid of $M$ generated by all idempotents, and further define
\begin{align*}
\F(M)&=\la E(M)\cup\G(M)\ra , & \GLR(M) &= \la \GL(M)\cup\GR(M)\ra, \\
\FL(M)&=\la E(M)\cup\GL(M)\ra , &\FLR(M) &= \la E(M)\cup\GL(M)\cup\GR(M)\ra \\
\FR(M)&=\la E(M)\cup\GR(M)\ra , & & = \la E(M)\cup\GLR(M)\ra.
\end{align*}
It will also be convenient to write 
\[
\I(M)=M \AND\O(M)=\{1_M\}.
\]

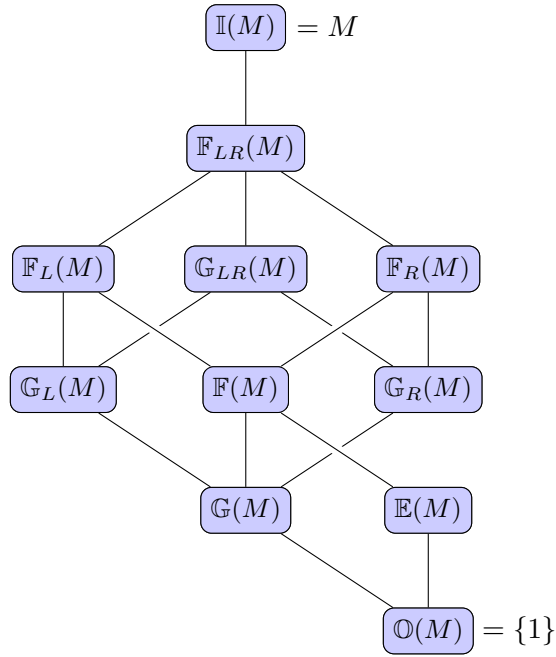
\begin{figure}[ht]
\begin{center}
\begin{tikzpicture}[scale=.8]
\node[rounded corners,rectangle,draw,fill=blue!20] (M) at (0,10) {\small $\I(M)$};
\node[rounded corners,rectangle,draw,fill=blue!20] (FLR) at (0,8) {\small $\F_{LR}(M)$};
\node[rounded corners,rectangle,draw,fill=blue!20] (FL) at (-3,6) {\small $\F_L(M)$};
\node[rounded corners,rectangle,draw,fill=blue!20] (GLR) at (0,6) {\small $\G_{LR}(M)$};
\node[rounded corners,rectangle,draw,fill=blue!20] (FR) at (3,6) {\small $\F_R(M)$};
\node[rounded corners,rectangle,draw,fill=blue!20] (GL) at (-3,4) {\small $\G_L(M)$};
\node[rounded corners,rectangle,draw,fill=blue!20] (F) at (0,4) {\small $\F(M)$};
\node[rounded corners,rectangle,draw,fill=blue!20] (GR) at (3,4) {\small $\G_R(M)$};
\node[rounded corners,rectangle,draw,fill=blue!20] (G) at (0,2) {\small $\G(M)$};
\node[rounded corners,rectangle,draw,fill=blue!20] (E) at (3,2) {\small $\E(M)$};
\node[rounded corners,rectangle,draw,fill=blue!20] (1) at (3,0) {\small $\O(M)$};
\draw (M)--(FLR)--(FL)--(GL)--(G)--(GR)--(FR)--(FLR);
\draw (GL)--(GLR)--(GR);
\draw[line width=1.5mm,white] (FL)--(F)--(FR);
\draw (FL)--(F)--(FR);
\draw (FLR)--(GLR);
\draw[line width=1.5mm,white] (F)--(E);
\draw (1)--(G)--(F)--(E)--(1);
\node () at (4.5,0) {$=\{1\}$};
\node () at (1.35,10.02) {$=M$};
\node[rounded corners,rectangle,draw,fill=blue!20] (M) at (0,10) {\small $\I(M)$};
\node[rounded corners,rectangle,draw,fill=blue!20] (FLR) at (0,8) {\small $\F_{LR}(M)$};
\node[rounded corners,rectangle,draw,fill=blue!20] (FL) at (-3,6) {\small $\F_L(M)$};
\node[rounded corners,rectangle,draw,fill=blue!20] (GLR) at (0,6) {\small $\G_{LR}(M)$};
\node[rounded corners,rectangle,draw,fill=blue!20] (FR) at (3,6) {\small $\F_R(M)$};
\node[rounded corners,rectangle,draw,fill=blue!20] (GL) at (-3,4) {\small $\G_L(M)$};
\node[rounded corners,rectangle,draw,fill=blue!20] (F) at (0,4) {\small $\F(M)$};
\node[rounded corners,rectangle,draw,fill=blue!20] (GR) at (3,4) {\small $\G_R(M)$};
\node[rounded corners,rectangle,draw,fill=blue!20] (G) at (0,2) {\small $\G(M)$};
\node[rounded corners,rectangle,draw,fill=blue!20] (E) at (3,2) {\small $\E(M)$};
\node[rounded corners,rectangle,draw,fill=blue!20] (1) at (3,0) {\small $\O(M)$};
\end{tikzpicture}
\caption{The generic shape of the lattice $\L(M)$.  In general these submonoids need not be distinct.}
\label{f:lattice}
\end{center}
\end{figure}

The relative containments of the submonoids defined above are shown in Figure \ref{f:lattice}.  Note that Figure~\ref{f:lattice} pictures the generic case, but that these submonoids need not be distinct in general; cf.~Figures \ref{f:L1}, \ref{f:L2} and \ref{f:L3}.

\subsection{Functors}\label{ss:fun}

We write $\M$ for the (locally small) category of all monoids.  The hom-set $\M(M,N)$ consists of all monoid homomorphisms $M\to N$ (each of which, recall, maps $1_M$ to $1_N$).

Now suppose $\X$ is one of $\E$, $\G$, $\GL$, $\GR$, $\GLR$, $\F$, $\FL$, $\FR$, $\FLR$, $\I$ or $\O$.  For any monoid $M$, $\X(M)$ is a submonoid of $M$, so it follows that $\X$ is an operator $\M\to\M$.  In fact, since any monoid homomorphism ${f:M\to N}$ maps idempotents (respectively, left units, right units, or units) of $M$ to idempotents (respectively, left units, right units, or units) of $N$, it is clear that~$f$ maps $\X(M)$ into $\X(N)$.  Thus, we may define~$\X(f):\X(M)\to\X(N)$ to be the restriction of $f$ to~$\X(M)$.  It then quickly follows that $\X$ is a functor~$\M\to\M$.  
We will write
\[
\Fu = \{\O,\E,\G,\GL,\GR,\GLR,\F,\FL,\FR,\FLR,\I\}
\]
for the set of all these functors.  

The direct product operation gives $\M$ the structure of a \emph{(symmetric) monoidal category}; see \cite[Chapters~VII and XI]{MacLane1998} and \cite{JS1993}.  The next lemma says that the functors from $\Fu$ are monoidal.  

\begin{lemma}\label{l:DP}
For any $\X\in\Fu$, and for any two monoids $M$ and $N$, we have
\[
\X(M\times N) = \X(M) \times \X(N).
\]
\end{lemma}

\pf
This is clear if $\X$ is $\O$ or $\I$.  For the other functors, it follows quickly from the fact that $(x,y)$ is an idempotent (or a left, right or two-sided unit) of $M\times N$ if and only if $x$ and $y$ are idempotents (or left, right or two-sided units) of $M$ and $N$, respectively.
\epf

It will also be convenient to record the following obvious fact.  For a monoid $M$, we write $M^0$ for the monoid obtained by adjoining a new zero element $0$ to $M$.  

\begin{lemma}\label{l:M0}
For any monoid $M$ we have
\[
\epfreseq
\X(M^0) = \begin{cases}
\X(M) &\text{if $\X$ is one of $\O$, $\G$, $\GL$, $\GR$ or $\GLR$}\\
\X(M)\cup\{0\} &\text{if $\X$ is one of $\I$, $\E$, $\F$, $\FL$, $\FR$ or $\FLR$.}
\end{cases}
\]
\end{lemma}

\subsection{Lattices}\label{ss:lat}

For a monoid $M$, we write 
\[
\L(M) = \set{\X(M)}{\X\in\Fu}
\]
for the set of all submonoids of $M$ defined in Section \ref{ss:sub}.  The set $\L(M)$ is partially ordered by inclusion; its Hasse diagram in the generic case is shown in Figure \ref{f:lattice}.  

We denote by $\Sub(M)$ the set of all submonoids of $M$, and we note that $\Sub(M)$ is a lattice with meet and join operations defined by
\[
S\wedge T=S\cap T \AND S\vee T=\la S\cup T\ra \qquad\text{for submonoids~$S$ and $T$ of $M$.}
\]
Throughout this article, the $\vee$ symbol will be used exclusively for the join operation in $\Sub(M)$.

\begin{prop}\label{p:L}
For any monoid $M$, the set $\L(M)$ is a finite $\vee$-subsemilattice of $\Sub(M)$, with top element~$\I(M)=M$ and bottom element~$\O(M)=\{1\}$.  Consequently, $\L(M)$ is a lattice.
\end{prop}

\pf
It is clear that $M$ and $\{1\}$ are the top and bottom elements of $\L(M)$.  Since a finite $\vee$-semilattice with a bottom element is a lattice (with the meet of two elements equal to the join of all common lower bounds), it suffices to show that $\L(M)$ is closed under $\vee$.  This is easily checked, using the definitions of the submonoids.  For example:
\[
\FL(M)\vee\GR(M) = \la \E(M)\cup\GL(M)\ra \vee \GR(M) = \la \E(M)\cup\GL(M)\cup\GR(M)\ra = \FLR(M).  \qedhere
\]
\epf

\begin{rem}\label{r:notlattice1}
The previous result did not say that $\L(M)$ is a sublattice of $\Sub(M)$ because this is not the case in general.  Specifically, $\L(M)$ is not always a $\wedge$-subsemilattice of $\Sub(M)$, meaning that the intersection of two submonoids from $\L(M)$ might not belong to~$\L(M)$; cf.~Remark \ref{r:notlattice2}.
\end{rem}

\subsection{Examples}\label{ss:egs}

Before we move on, we pause to consider some basic examples.  These should serve to illustrate the above ideas, but will also be useful later for proving some of our main results.

First, if $G$ is a group, then clearly every element is a (two-sided) unit, and the only idempotent is the identity element.  It quickly follows that the submonoids $\X(G)$, $\X\in\Fu$, are as listed in the first column of Table \ref{t:GNB_small}.

Next, suppose $E$ is an idempotent-generated monoid.  Clearly $\E(E)=E$.  It follows from \cite[Lemma~2.1]{JE_IBM} that $\GL(E)=\GR(E)=\G(E)=\{1\}$.  Thus, the submonoids $\X(G)$, $\X\in\Fu$, are as listed in the second column of Table \ref{t:GNB_small}.

Next, we denote by $\Pos=\{1,2,3,\ldots\}$ the multiplicative monoid of positive integers.  This time, $1$ is the unique unit, and also the unique idempotent.  The submonoids $\X(\Pos)$, $\X\in\Fu$, are listed in the third column of Table \ref{t:GNB_small}.

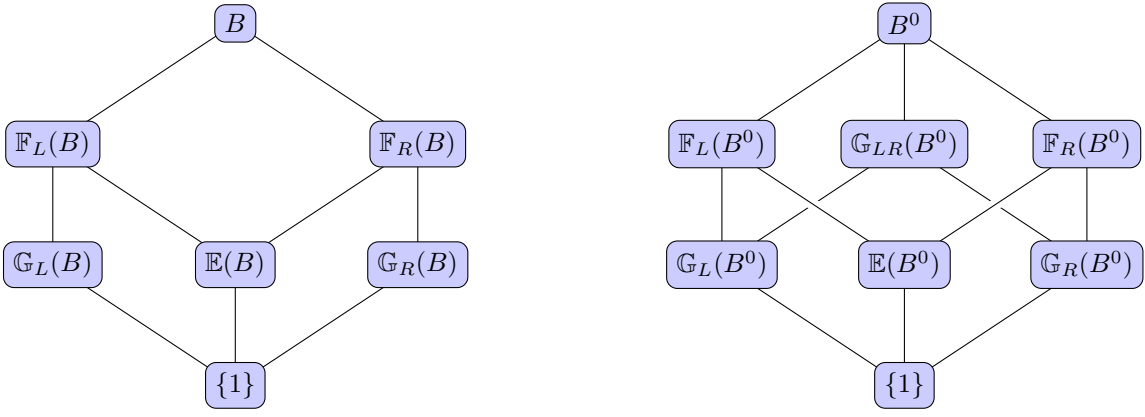
\begin{figure}[ht]
\begin{center}
\begin{tikzpicture}[scale=.8]
\begin{scope}
\node[rounded corners,rectangle,draw,fill=blue!20] (FLR) at (0,8) {};
\node[rounded corners,rectangle,draw,fill=blue!20] (FL) at (-3,6) {};
\node[rounded corners,rectangle,draw,fill=blue!20] (FR) at (3,6) {};
\node[rounded corners,rectangle,draw,fill=blue!20] (GL) at (-3,4) {};
\node[rounded corners,rectangle,draw,fill=blue!20] (F) at (0,4) {};
\node[rounded corners,rectangle,draw,fill=blue!20] (GR) at (3,4) {};
\node[rounded corners,rectangle,draw,fill=blue!20] (G) at (0,2) {};
\draw (FLR)--(FL)--(GL)--(G)--(GR)--(FR)--(FLR);
\draw (FL)--(F)--(FR) (G)--(F);
\node[rounded corners,rectangle,draw,fill=blue!20] (FLR) at (0,8) {\small $B$};
\node[rounded corners,rectangle,draw,fill=blue!20] (FL) at (-3,6) {\small $\FL(B)$};
\node[rounded corners,rectangle,draw,fill=blue!20] (FR) at (3,6) {\small $\FR(B)$};
\node[rounded corners,rectangle,draw,fill=blue!20] (GL) at (-3,4) {\small $\GL(B)$};
\node[rounded corners,rectangle,draw,fill=blue!20] (F) at (0,4) {\small $\E(B)$};
\node[rounded corners,rectangle,draw,fill=blue!20] (GR) at (3,4) {\small $\GR(B)$};
\node[rounded corners,rectangle,draw,fill=blue!20] (G) at (0,2) {\small $\{1\}$};
\end{scope}
\begin{scope}[shift={(11,0)}]
\node[rounded corners,rectangle,draw,fill=blue!20] (FLR) at (0,8) {};
\node[rounded corners,rectangle,draw,fill=blue!20] (FL) at (-3,6) {};
\node[rounded corners,rectangle,draw,fill=blue!20] (GLR) at (0,6) {};
\node[rounded corners,rectangle,draw,fill=blue!20] (FR) at (3,6) {};
\node[rounded corners,rectangle,draw,fill=blue!20] (GL) at (-3,4) {};
\node[rounded corners,rectangle,draw,fill=blue!20] (F) at (0,4) {};
\node[rounded corners,rectangle,draw,fill=blue!20] (GR) at (3,4) {};
\node[rounded corners,rectangle,draw,fill=blue!20] (G) at (0,2) {};
\draw (FLR)--(FL)--(GL)--(G)--(GR)--(FR)--(FLR);
\draw (GL)--(GLR)--(GR);
\draw[line width=1.5mm,white] (FL)--(F)--(FR) (G)--(F);
\draw (FL)--(F)--(FR) (G)--(F);
\draw (FLR)--(GLR);
\node[rounded corners,rectangle,draw,fill=blue!20] (FLR) at (0,8) {\small $B^0$};
\node[rounded corners,rectangle,draw,fill=blue!20] (FL) at (-3,6) {\small $\FL(B^0)$};
\node[rounded corners,rectangle,draw,fill=blue!20] (GLR) at (0,6) {\small $\GLR(B^0)$};
\node[rounded corners,rectangle,draw,fill=blue!20] (FR) at (3,6) {\small $\FR(B^0)$};
\node[rounded corners,rectangle,draw,fill=blue!20] (GL) at (-3,4) {\small $\GL(B^0)$};
\node[rounded corners,rectangle,draw,fill=blue!20] (F) at (0,4) {\small $\E(B^0)$};
\node[rounded corners,rectangle,draw,fill=blue!20] (GR) at (3,4) {\small $\GR(B^0)$};
\node[rounded corners,rectangle,draw,fill=blue!20] (G) at (0,2) {\small $\{1\}$};
\end{scope}
\end{tikzpicture}
\caption{The lattices $\L(B)$ and $\L(B^0)$, where $B$ is the bicyclic monoid.  In both diagrams, the nodes represent distinct submonoids.}
\label{f:B}
\end{center}
\end{figure}

The \emph{bicyclic monoid} $B$ is defined by the monoid presentation $B=\la a,b:ba=1\ra$.  Because of the relation $ba=1$, we may think of the elements of $B$ as words of the form $a^mb^n$, where $m,n\geq0$.  Two such words $a^mb^n$ and $a^kb^l$ represent the same element of $B$ if and only if $m=k$ and $n=l$, and the product in $B$ is given by
\begin{equation}\label{e:B}
a^mb^n\cdot a^kb^l = a^{m+\mu-n}b^{l+\mu-k} \qquad\text{where $\mu=\max(n,k)$.}
\end{equation}
Any monoid generated by two elements $x,y$ for which $yx=1\not=xy$ is isomorphic to~$B$; see \cite[pp.~31--32]{Howie} for more details.  Idempotents of $B$ are words of the form $a^mb^m$ ($m\geq0$), and it is easily checked that idempotents commute, so that $\E(B)=E(B)$.  
Using \eqref{e:B}, it is easy to see that
\[
a^mb^n\cdot a^kb^l=1 \iff m=l=0 \text{ and } n=k,
\]
so that
\[
\GL(B)=\la a\ra=\{1,a,a^2,\ldots\} \AND \GR(B)=\la b\ra=\{1,b,b^2,\ldots\}.
\]
The fourth column of Table \ref{t:GNB_small} lists the submonoids $\X(B)$, $\X\in\Fu$; verification for the submonoids not discussed so far is an exercise.  
The fifth column of Table \ref{t:GNB_small} lists the corresponding submonoids of $B^0$ (the bicyclic monoid with a zero adjoined); cf.~Lemma \ref{l:M0}.
The lattices $\L(B)$ and $\L(B^0)$ are pictured in Figure~\ref{f:B}.

\begin{table}[ht]
\begin{center}
\begin{tabular}{l|ccccc}
$\X$ & $\X(G)$ & $\X(E)$ & $\X(\Pos)$ & $\X(B)$ & $\X(B^0)$ \\
\hline
$\O$ & $\{1\}$ & $\{1\}$ & $\{1\}$ & $\{1\}$  & $\{1\}$ \\
$\E$ & $\{1\}$ & $E$ & $\{1\}$ & $\set{a^mb^m}{m\geq0}$ & $\set{a^mb^m}{m\geq0}\cup\{0\}$ \\
$\G$ & $G$ & $\{1\}$ & $\{1\}$ & $\{1\}$ & $\{1\}$ \\
$\GL$ & $G$ & $\{1\}$ & $\{1\}$ & $\la a\ra$ & $\la a\ra$ \\
$\GR$ & $G$ & $\{1\}$ & $\{1\}$ & $\la b\ra$ & $\la b\ra$ \\
$\GLR$ & $G$ & $\{1\}$ & $\{1\}$ & $B$ & $B$ \\
$\F$ & $G$ & $E$ & $\{1\}$ & $\set{a^mb^m}{m\geq0}$ & $\set{a^mb^m}{m\geq0}\cup\{0\}$ \\
$\FL$ & $G$ & $E$ & $\{1\}$ & $\set{a^mb^n}{m\geq n}$ & $\set{a^mb^n}{m\geq n}\cup\{0\}$ \\
$\FR$ & $G$ & $E$ & $\{1\}$ & $\set{a^mb^n}{m\leq n}$ & $\set{a^mb^n}{m\leq n}\cup\{0\}$ \\
$\FLR$ & $G$ & $E$ & $\{1\}$  & $B$ &$B^0$ \\
$\I$ & $G$ & $E$ & $\Pos$ & $B$ &$B^0$ \\
\end{tabular}
\caption{The submonoids $\X(M)$, $\X\in\Fu$, for $M=G$ (a group), $M=E$ (an idempotent-generated monoid), $M=\Pos$ (the positive integers under multiplication), $M=B$ (the bicyclic monoid) and $M=B^0$ (the bicyclic monoid with a zero adjoined).}
\label{t:GNB_small}
\end{center}
\end{table}

\begin{rem}\label{r:notlattice2}
Proposition \ref{p:L} showed that the lattice $\L(M)$ is a $\vee$-subsemilattice of $\Sub(M)$, and we claimed in Remark \ref{r:notlattice1} that $\L(M)$ is not always a $\wedge$-subsemilattice.  We can use the above example of $M=B^0$ to verify this.  Indeed, using Table \ref{t:GNB_small} we see that the meet in $\Sub(B^0)$ of the submonoids $\FL(B^0)$ and $\GLR(B^0)$ is
\[
\FL(B^0) \cap \GLR(B^0) = \set{a^mb^n}{m\geq n},
\]
which does not belong to $\L(B^0)$.  Of course, the submonoids $\FL(B^0)$ and $\GLR(B^0)$ do have a meet in $\L(B^0)$ itself, as the latter is a lattice, but this meet in $\L(B^0)$ is $\GL(B^0)=\la a\ra$; cf.~Figure \ref{f:B}.

On the other hand, the lattice $\L(B)$ is a sublattice of $\Sub(B)$, as may be easily verified using Table \ref{t:GNB_small}.
\end{rem}

\section{Preliminary results}\label{s:prelim}

We now gather a number of technical results that will be useful in subsequent sections.  Section~\ref{ss:int} concerns intersections of various submonoids from the lattice $\L(M)$, and Section \ref{ss:collapse} concerns equalities between such submonoids.  Section \ref{ss:Green} establishes connections with Green's relations, in particular with stability (or otherwise) of the identity element.

Throughout this section, unless otherwise stated, $M$ will denote an arbitrary monoid.  It will also be convenient to abbreviate the submonoids $\X(M)$, $\X\in\Fu$, in obvious ways.  Specifically, we will often write
\begin{align}
\label{e:G} && G&=\G(M), & G_L&=\GL(M), & G_R&=\GR(M), & G_{LR}&=\GLR(M), \\
\label{e:F} E&=\E(M), & F&=\F(M), & F_L&=\FL(M), & F_R&=\FR(M), & F_{LR}&=\FLR(M).
\end{align}
A further piece of notation will also be convenient.  Often we will wish to give a statement or argument that holds regardless of subscripts, so will sometimes write $\G_\heartsuit(M)$ or $\G_\diamondsuit(M)$ to stand for any of $\G(M)$, $\GL(M)$, $\GR(M)$ or $\GLR(M)$.  Similarly, we will at times write $G_\heartsuit$ or $F_\diamondsuit$, etc.

\subsection{Intersections}\label{ss:int}

The next two results concern intersections of various submonoids of $M$.  We will sometimes make use of them without explicit reference.  The first concerns intersections with $E=\E(M)$.

\begin{lemma}\label{l:GE}
For any monoid $M$ we have
\ben
\item \label{GE1} $E\cap G=E\cap G_L=E\cap G_R=\{1\}$,
\item \label{GE2} $E\cap F=E\cap F_L=E\cap F_R=E\cap F_{LR}=E$.
\een
\end{lemma}

\pf
Part \ref{GE1} is part of \cite[Lemma 2.1]{JE_IBM}.  Part \ref{GE2} is clear, since $E$ is contained in each of $F,F_L,F_R,F_{LR}$.
\epf

\begin{rem}
The previous result did not say anything about $E\cap G_{LR}$.  Certainly $\{1\}\sub E\cap G_{LR} \sub E$, but we cannot say any more than this in general, since any of the following situations are possible (cf.~Table \ref{t:GNB_small}):
\bit
\item
$\{1\}= E\cap G_{LR} = E$: e.g., if $M$ is a group,
\item
$\{1\}= E\cap G_{LR} \subsetneq E$: e.g., if $M$ is a group with a zero adjoined,
\item
$\{1\}\subsetneq E\cap G_{LR} = E$: e.g., if $M$ is the bicyclic monoid,
\item
$\{1\}\subsetneq E\cap G_{LR} \subsetneq E$: e.g., if $M$ is the bicyclic monoid with a zero adjoined.
\eit
\end{rem}

The next result concerns intersections with $G_L=\GL(M)$.  There is an obvious dual result concerning intersections with $G_R=\GR(M)$, but we will not state it.

\begin{lemma}\label{l:GL}
For any monoid $M$ we have
\ben
\item \label{GL1} $G_L\cap G_L = G_L\cap G_{LR} = G_L\cap F_L = G_L\cap F_{LR} = G_L$,
\item \label{GL2} $G_L\cap G = G_L\cap G_R = G_L\cap F = G_L\cap F_R = G$.
\een
\end{lemma}

\pf
\firstpfitem{\ref{GL1}}  This is clear, since $G_L$ is contained in each of $G_L,G_{LR},F_L,F_{LR}$.

\pfitem{\ref{GL2}}  Since each of the stated intersections contains $G$, and since each of $G,G_R,F$ are contained in~$F_R$, it suffices to show that $G_L\cap F_R\sub G$.  To do so, suppose $x\in G_L\cap F_R$.  Since $x\in G_L$ we have $1=ax$ for some $a\in M$.  Since $F_R=G_R E$ by \cite[Lemma 2.5]{JE_IBM}, we also have $x=ge_1\cdots e_k$ for some $k\geq0$, and some $g\in G_R$ and $e_1,\ldots,e_k\in E(M)$.  We may assume that $k$ is minimal among all such expressions; in particular, $e_i\not=1$ for all $1\leq i\leq k$ .  If $k\geq1$, then $e_k\not=1$ and $x=xe_k$, which gives $e_k = 1e_k=axe_k=ax=1$, a contradiction.  Thus, $k=0$, so that $x=g\in G_R$.  It follows that $x\in G_L\cap G_R=G$, as required.
\epf

\subsection{Collapse}\label{ss:collapse}

We have already observed that the submonoids of $M$ defined in Section \ref{ss:sub} are not always distinct.  Roughly speaking, this means that certain ``collapse'' can occur in the lattice $\L(M)$.  The next two results show that such collapse happens in a somewhat controlled manner, in the sense that equalities between certain submonoids imply other such equalities.

\begin{lemma}\label{l:GLR}
For a monoid $M$, the following are equivalent:
\ben
\item \label{GLR1} $G$, $G_L$, $G_R$ and $G_{LR}$ are not all equal,
\item \label{GLR2} $G$, $G_L$, $G_R$ and $G_{LR}$ are pairwise distinct,
\item \label{FLR1} $F$, $F_L$, $F_R$ and $F_{LR}$ are not all equal,
\item \label{FLR2} $F$, $F_L$, $F_R$ and $F_{LR}$ are pairwise distinct,
\item \label{GLR3} $G_{LR}$ contains infinitely many idempotents,
\item \label{GLR4} $G_{LR}$ contains a nontrivial idempotent,
\item \label{GLR5} $E\cap G_{LR}\not=\{1\}$.
\een
\end{lemma}

\pf
We begin by establishing the equivalence of items involving submonoids of the form $G_\hs$.

\pfitem{\ref{GLR1}$\implies$\ref{GLR2}}  We prove the contrapositive: i.e., that if any two of the stated submonoids are equal, then all four are equal.  
\bit
\item If $G=G_L$, then also $G=G_R$ (cf.~\cite[Lemma 2.3]{JE_IBM}), and $G_{LR}= G_L\vee G_R= G\vee G=G$.  
\item The $G=G_R$ case is dual.
\item If $G=G_{LR}$, then $G=G_L\cap G=G_L\cap G_{LR}=G_L$, reducing to the first case.
\item If $G_L=G_R$, then $G=G_L\cap G_R=G_L\cap G_L=G_L$.
\item If $G_L=G_{LR}$, then $G=G_L\cap G_R=G_{LR}\cap G_R=G_R$.
\item The $G_R=G_{LR}$ case is again dual.
\eit

\pfitem{\ref{GLR2}$\implies$\ref{GLR3}}  Suppose $G\not=G_L$, and let $x\in G_L\sm G$ be arbitrary.  Then $1=ax$ for some $a\in M$, and we note that $a\in G_R$.  Since $x\not\in G$ we have $xa\not=1$.  It follows that $\la a,x\ra$ is bicyclic, and hence contains infinitely many idempotents (of the form $x^ma^m$ for each $m\geq0$).  Since $x\in G_L$ and $a\in G_R$, it follows that $\la a,x\ra\sub G_{LR}$.

\pfitem{\ref{GLR3}$\implies$\ref{GLR4} and \ref{GLR4}$\implies$\ref{GLR5}}  These are clear.

\pfitem{\ref{GLR5}$\implies$\ref{GLR1}}  If $E\cap G_{LR}\not=\{1\}$, then $G_{LR}\not=G$ because $E\cap G=\{1\}$.

\afterpfitem Now that we know \ref{GLR1}, \ref{GLR2}, \ref{GLR3}--\ref{GLR5} are equivalent, it is time to tie these in with \ref{FLR1} and \ref{FLR2}.

\pfitem{\ref{GLR2}$\implies$\ref{FLR2}}  Suppose \ref{GLR2} holds.  From Lemma \ref{l:GL}, we have
\[
G_L\cap F_L = G_L\cap F_{LR} = G_L \not= G = G_L\cap F = G_L\cap F_R,
\]
and it follows that $\{F_L,F_{LR}\} \cap \{F,F_R\} = \emptyset$.  We similarly obtain $\{F_R,F_{LR}\} \cap \{F,F_L\} = \emptyset$ from the dual of Lemma \ref{l:GL}.

\pfitem{\ref{FLR2}$\implies$\ref{FLR1}}  This is clear.

\pfitem{\ref{FLR1}$\implies$\ref{GLR1}}  Aiming to prove the contrapositive, suppose $G_\hs=G_\ds$ for distinct subscripts $\hs,\ds$.  Then $F_\hs= E\vee G_\hs= E\vee G_\ds = F_\ds$.
\epf

The previous lemma concerned collapse in $\L(M)$ within the two ``diamonds'' $\{G,G_L,G_R,G_{LR}\}$ and $\{F,F_L,F_R,F_{LR}\}$.  The next concerns collapse at the very bottom of the lattice, namely between $\{1\}$ and $G$ or $E$.  In particular, it shows that $E=\{1\}$ has the significant consequence of collapsing the whole ``cube'' section of the lattice: i.e., the interval from~$G$ to $F_{LR}$.

\begin{lemma}\label{l:EG}
For any monoid $M$ we have
\ben
\item \label{EG2} $G=\{1\} \iff F=E$,
\item \label{EG1} $E=\{1\} \iff F_{LR}=G \iff \{G,G_L,G_R\} \cap \{F,F_L,F_R,F_{LR}\} \not=\emptyset$.
\een
\end{lemma}

\pf
\firstpfitem{\ref{EG2}}  If $G=\{1\}$ then $F= E\vee G= E\vee\{1\}=E$.  Conversely, if $F=E$, then since $G\sub F$, we have $G=F\cap G=E\cap G=\{1\}$.

\pfitem{\ref{EG1}}  If $E=\{1\}$, then $G_{LR}$ contains no nontrivial idempotents, so by Lemma \ref{l:GLR} we have $G=G_{LR}$; but then $F_{LR}= E\vee G_{LR} = \{1\}\vee G=G$.

If $F_{LR}=G$, then obviously $\{G,G_L,G_R\} \cap \{F,F_L,F_R,F_{LR}\} \not=\emptyset$.

Finally, suppose the two stated sets of submonoids have nonempty intersection, say $G_\hs=F_\ds$, noting that $\hs\not=LR$.  Then Lemma \ref{l:GE} gives $E=E\cap F_\ds=E\cap G_\hs=\{1\}$.
\epf

\begin{rem}
The submonoid $G_{LR}=\GLR(M)$ was not mentioned in Lemma \ref{l:EG}\ref{EG1}, since it is possible to have $G_{LR}=F_{LR}$ but $E\not=\{1\}$.  For example, this happens when $M$ is the bicyclic monoid; cf.~Table \ref{t:GNB_small} and Figure \ref{f:B}.
\end{rem}

\subsection{Green's relations and stability}\label{ss:Green}

Recall that for elements $x$ and $y$ of a monoid $M$, we write
\[
x\gL y\iff Mx=My \COMMA
x\gR y\iff xM=yM \COMMA
x\gJ y\iff MxM=MyM .
\]
We also set ${\gH}={\gL}\cap{\gR}$ and ${\gD}={\gL}\vee{\gR}$ (the join in the lattice of equivalences).  
These five equivalences,~$\gL$,~$\gR$,~$\gJ$,~$\gH$ and~$\gD$, are called \emph{Green's relations} \cite{Green1951}, and are essential tools in semigroup theory.  
Equivalent formulations in terms of divisibility may also be given; for example, $x\gL y$ if and only if $x=ay$ and $y=bx$ for some $a,b\in M$.  See \cite[Chapter~2]{CPbook} or \cite[Chapter~2]{Howie} for more background on Green's relations.

If $\gK$ is one of Green's relations, we denote by~$K_x=\set{y\in M}{x\gK y}$ the $\gK$-class of $x\in M$.  One may easily check that the submonoids consisting of one- or two-sided units are certain Green's classes containing the identity:
\[
G_L=\GL(M) = L_1 \COMMA G_R=\GR(M) = R_1 \COMMA G=\G(M) = H_1.
\]

An element $x$ of a monoid $M$ is \emph{stable} if the following implications hold for all $a\in M$:
\begin{equation}\label{e:stab}
ax\gJ x \implies ax\gL x \AND xa\gJ x \implies xa\gR x.
\end{equation}
If $x$ is not stable, we will call it \emph{unstable}.  For more on stability, see \cite[Section 2.3]{Lallement1979}, \cite[Section A.2]{RSbook} or \cite{EH2019}.

Taking $x=1$ to be the identity element in \eqref{e:stab}, and keeping in mind that ${\gH}={\gL}\cap{\gR}$, we see that $1$ is stable if and only if $a\gJ1\implies a\gH1$ for all $a\in M$.  This implication is equivalent to $J_1\sub H_1$.  Since $H_x\sub J_x$ for any $x$, it follows that $1$ is stable if and only if $J_1=H_1$: i.e., $J_1=G$.

\newpage

\begin{lemma}\label{l:J}
For a monoid $M$, the following are equivalent:
\ben
\item \label{J1} $G=G_L$,
\item \label{J2} $G=G_R$,
\item \label{J3} $G=J_1$,
\item \label{J4} $H_1=L_1=R_1=D_1=J_1$,
\item \label{J5} $M$ has no bicyclic submonoid,
\item \label{J6} the identity element $1$ is stable.
\een
\end{lemma}

\pf
\firstpfitem{\ref{J1}$\iff$\ref{J2}$\iff$\ref{J5}}  These are part of \cite[Lemma 2.3]{JE_IBM}; cf.~Lemma \ref{l:GLR}.  

\pfitem{\ref{J3}$\iff$\ref{J4}}  This follows from $G=H_1\sub L_1,R_1\sub D_1\sub J_1$, which itself follows from ${\gH}\sub{\gL},{\gR}\sub{\gD}\sub{\gJ}$.

\pfitem{\ref{J3}$\iff$\ref{J6}}  This was discussed before the statement of the lemma.

\pfitem{\ref{J1}$\implies$\ref{J3}}  If $G=G_L$ holds, then so too does $G=G_R$ (as \ref{J1}$\iff$\ref{J2}).  Since $G=H_1\sub J_1$, it is enough to show that $J_1\sub G$.  To do so, let $x\in J_1$.  Then $1=axb$ for some $a,b\in M$.  Since $1=a(xb)$ we have $a\in G_R=G$, and similarly $b\in G$.  But then $x=a^{-1}(axb)b^{-1}=a^{-1}b^{-1}\in G$.

\pfitem{\ref{J4}$\implies$\ref{J1}}  If \ref{J4} holds, then $G=H_1=L_1=G_L$.
\epf

\begin{rem}\label{r:dichotomy}
The second condition of Lemma \ref{l:GLR} and the first condition of Lemma \ref{l:J} are clearly mutually exclusive.  It follows that a monoid either satisfies all of the conditions of Lemma \ref{l:GLR} and none of the conditions of Lemma \ref{l:J}, or vice versa.  This yields a dichotomy that will allow a convenient split in the argument of the next section.
\end{rem}

\section{Classification of lattice invariants}\label{s:class}

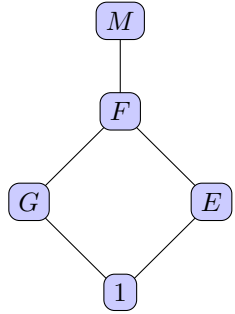
\begin{figure}[ht]
\begin{center}
\begin{tikzpicture}[scale=.6]
\node[rounded corners,rectangle,draw,fill=blue!20] (M) at (0,6) {};
\node[rounded corners,rectangle,draw,fill=blue!20] (F) at (0,4) {};
\node[rounded corners,rectangle,draw,fill=blue!20] (G) at (-2,2) {};
\node[rounded corners,rectangle,draw,fill=blue!20] (E) at (2,2) {};
\node[rounded corners,rectangle,draw,fill=blue!20] (1) at (0,0) {};
\draw (F)--(G)--(1)--(E)--(F)--(M);
\node[rounded corners,rectangle,draw,fill=blue!20] (M) at (0,6) {\small $M$};
\node[rounded corners,rectangle,draw,fill=blue!20] (F) at (0,4) {\small $F$};
\node[rounded corners,rectangle,draw,fill=blue!20] (G) at (-2,2) {\small $G$};
\node[rounded corners,rectangle,draw,fill=blue!20] (E) at (2,2) {\small $E$};
\node[rounded corners,rectangle,draw,fill=blue!20] (1) at (0,0) {\small $1$};
\end{tikzpicture}
\caption{The generic shape of the lattice $\L(M)$ when $M$ has a stable identity.}
\label{f:lattice_small}
\end{center}
\end{figure}

\begin{figure}[ht]
\begin{center}
\begin{tikzpicture}[scale=.5]
\begin{scope}[shift={(0,0)}]
\node[rounded corners,rectangle,draw,fill=blue!20] (M) at (0,6) {};
\node[rounded corners,rectangle,draw,fill=blue!20] (F) at (0,4) {};
\node[rounded corners,rectangle,draw,fill=blue!20] (G) at (-2,2) {};
\node[rounded corners,rectangle,draw,fill=blue!20] (E) at (2,2) {};
\node[rounded corners,rectangle,draw,fill=blue!20] (1) at (0,0) {};
\draw (F)--(G)--(1)--(E)--(F)--(M);
\node[rounded corners,rectangle,draw,fill=blue!20] (M) at (0,6) {\small $M$};
\node[rounded corners,rectangle,draw,fill=blue!20] (F) at (0,4) {\small $F$};
\node[rounded corners,rectangle,draw,fill=blue!20] (G) at (-2,2) {\small $G$};
\node[rounded corners,rectangle,draw,fill=blue!20] (E) at (2,2) {\small $E$};
\node[rounded corners,rectangle,draw,fill=blue!20] (1) at (0,0) {\small $1$};
\node () at (0,7.2) {$(1,0,0,0)$};
\end{scope}
\begin{scope}[shift={(8,0)}]
\node[rounded corners,rectangle,draw,fill=blue!20] (M) at (0,6) {\small $M$};
\node[rounded corners,rectangle,draw,fill=blue!20] (F) at (0,4) {\small $F=E$};
\node[rounded corners,rectangle,draw,fill=blue!20] (G) at (0,2) {\small $G=1$};
\draw (G)--(F)--(M);
\node () at (0,7.2) {$(1,0,0,1)$};
\end{scope}
\begin{scope}[shift={(16,0)}]
\node[rounded corners,rectangle,draw,fill=blue!20] (M) at (0,6) {\small $M$};
\node[rounded corners,rectangle,draw,fill=blue!20] (F) at (0,4) {\small $F=G$};
\node[rounded corners,rectangle,draw,fill=blue!20] (E) at (0,2) {\small $E=1$};
\draw (E)--(F)--(M);
\node () at (0,7.2) {$(1,0,1,0)$};
\end{scope}
\begin{scope}[shift={(24,0)}]
\node[rounded corners,rectangle,draw,fill=blue!20] (M) at (0,6) {\small $M$};
\node[rounded corners,rectangle,draw,fill=blue!20] (F) at (0,4) {\small $F=G=E=1$};
\draw (F)--(M);
\node () at (0,7.2) {$(1,0,1,1)$};
\end{scope}
\begin{scope}[shift={(0,-8)}]
\node[rounded corners,rectangle,draw,fill=blue!20] (F) at (0,4) {};
\node[rounded corners,rectangle,draw,fill=blue!20] (G) at (-2,2) {};
\node[rounded corners,rectangle,draw,fill=blue!20] (E) at (2,2) {};
\node[rounded corners,rectangle,draw,fill=blue!20] (1) at (0,0) {};
\draw (F)--(G)--(1)--(E)--(F);
\node[rounded corners,rectangle,draw,fill=blue!20] (F) at (0,4) {\small $M=F$};
\node[rounded corners,rectangle,draw,fill=blue!20] (G) at (-2,2) {\small $G$};
\node[rounded corners,rectangle,draw,fill=blue!20] (E) at (2,2) {\small $E$};
\node[rounded corners,rectangle,draw,fill=blue!20] (1) at (0,0) {\small $1$};
\node () at (0,5.2) {$(1,1,0,0)$};
\end{scope}
\begin{scope}[shift={(8,-8)}]
\node[rounded corners,rectangle,draw,fill=blue!20] (F) at (0,4) {\small $M=F=E$};
\node[rounded corners,rectangle,draw,fill=blue!20] (G) at (0,2) {\small $G=1$};
\draw (G)--(F);
\node () at (0,5.2) {$(1,1,0,1)$};
\end{scope}
\begin{scope}[shift={(16,-8)}]
\node[rounded corners,rectangle,draw,fill=blue!20] (F) at (0,4) {\small $M=F=G$};
\node[rounded corners,rectangle,draw,fill=blue!20] (E) at (0,2) {\small $E=1$};
\draw (E)--(F);
\node () at (0,5.2) {$(1,1,1,0)$};
\end{scope}
\begin{scope}[shift={(24,-8)}]
\node[rounded corners,rectangle,draw,fill=blue!20] (F) at (0,4) {\small $M=F=G=E=1$};
\node () at (0,5.2) {$(1,1,1,1)$};
\end{scope}
\end{tikzpicture}
\caption{The lattice $\L(M)$ when $M$ has a stable identity, according to the type $\T(M)=(1,i,j,k)$.  In each case, the nodes represent distinct submonoids of $M$.}
\label{f:L1}
\end{center}
\end{figure}
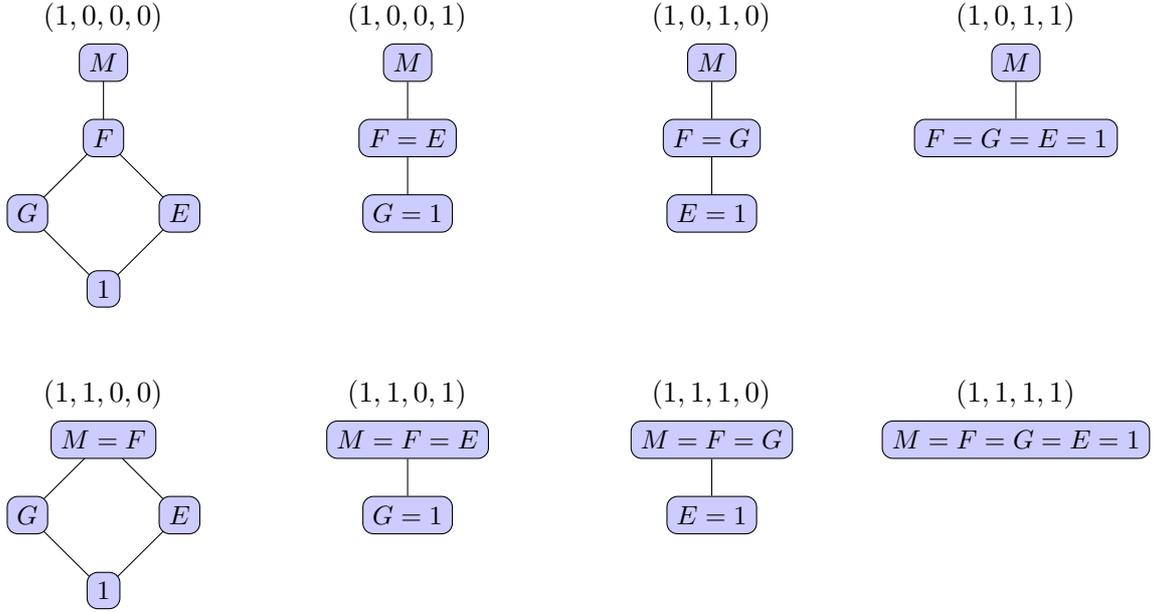

\begin{figure}[ht]
\begin{center}
\begin{tikzpicture}[scale=.6]
\node (M) at (0,10) {};
\node (FLR) at (0,8) {};
\node (FL) at (-3,6) {};
\node (GLR) at (0,6){};
\node (FR) at (3,6) {};
\node (GL) at (-3,4) {};
\node (F) at (0,4) {};
\node (GR) at (3,4) {};
\node (G) at (0,2) {};
\node (E) at (3,2) {};
\node (1) at (3,0) {};
\draw[ultra thick] (FLR)--(FL)--(GL)--(G)--(GR)--(FR)--(FLR);
\draw[ultra thick] (GL)--(GLR)--(GR);
\draw[line width=2mm,white] (FL)--(F)--(FR);
\draw[ultra thick] (FL)--(F)--(FR);
\draw[ultra thick] (G)--(F) (E)--(1);
\draw (M)--(FLR)--(GLR);
\draw[line width=1.5mm,white] (F)--(E);
\draw (F)--(E) (G)--(1);
\node[rounded corners,rectangle,draw,fill=blue!20] (M) at (0,10) {\small $M$};
\node[rounded corners,rectangle,draw,fill=red!20] (FLR) at (0,8) {\small $F_{LR}$};
\node[rounded corners,rectangle,draw,fill=red!20] (FL) at (-3,6) {\small $F_L$};
\node[rounded corners,rectangle,draw,fill=blue!20] (GLR) at (0,6) {\small $G_{LR}$};
\node[rounded corners,rectangle,draw,fill=red!20] (FR) at (3,6) {\small $F_R$};
\node[rounded corners,rectangle,draw,fill=red!20] (GL) at (-3,4) {\small $G_L$};
\node[rounded corners,rectangle,draw,fill=red!20] (F) at (0,4) {\small $F$};
\node[rounded corners,rectangle,draw,fill=red!20] (GR) at (3,4) {\small $G_R$};
\node[rounded corners,rectangle,draw,fill=red!20] (G) at (0,2) {\small $G$};
\node[rounded corners,rectangle,draw,fill=blue!20] (E) at (3,2) {\small $E$};
\node[rounded corners,rectangle,draw,fill=blue!20] (1) at (3,0) {\small $1$};
\end{tikzpicture}
\caption{The generic shape of the lattice $\L(M)$ when $M$ has an unstable identity.  The submonoids shaded red are distinct, and thick lines indicate proper containment.}
\label{f:lattice_red}
\end{center}
\end{figure}
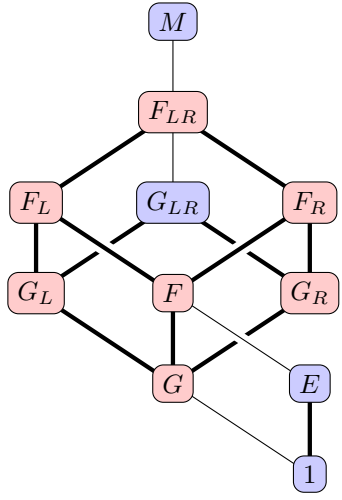

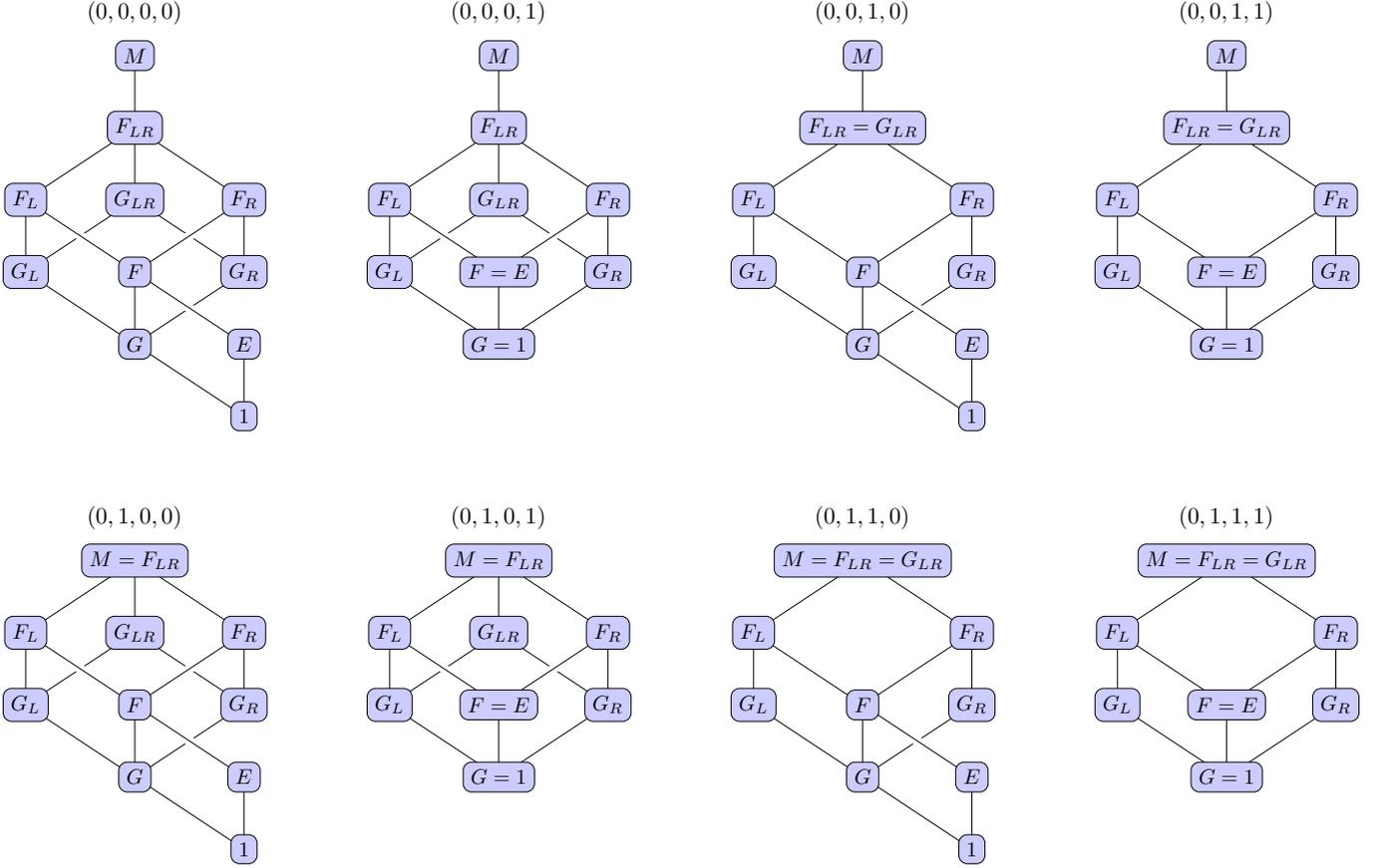
\begin{figure}[ht]
\begin{center}
\scalebox{0.8}{
\begin{tikzpicture}[scale=.6]
\begin{scope}[shift={(0,0)}]
\node[rounded corners,rectangle,draw,fill=blue!20] (M) at (0,10) {};
\node[rounded corners,rectangle,draw,fill=blue!20] (FLR) at (0,8) {};
\node[rounded corners,rectangle,draw,fill=blue!20] (FL) at (-3,6) {};
\node[rounded corners,rectangle,draw,fill=blue!20] (GLR) at (0,6) {};
\node[rounded corners,rectangle,draw,fill=blue!20] (FR) at (3,6) {};
\node[rounded corners,rectangle,draw,fill=blue!20] (GL) at (-3,4) {};
\node[rounded corners,rectangle,draw,fill=blue!20] (F) at (0,4) {};
\node[rounded corners,rectangle,draw,fill=blue!20] (GR) at (3,4) {};
\node[rounded corners,rectangle,draw,fill=blue!20] (G) at (0,2) {};
\node[rounded corners,rectangle,draw,fill=blue!20] (E) at (3,2) {};
\node[rounded corners,rectangle,draw,fill=blue!20] (1) at (3,0) {};
\draw (FLR)--(FL)--(GL)--(G)--(GR)--(FR)--(FLR);
\draw (GL)--(GLR)--(GR);
\draw[line width=1.5mm,white] (FL)--(F)--(FR);
\draw (FL)--(F)--(FR);
\draw (G)--(F) (E)--(1);
\draw (M)--(FLR)--(GLR);
\draw[line width=1.5mm,white] (F)--(E);
\draw (F)--(E) (G)--(1);
\node () at (0,11.2) {$(0,0,0,0)$};
\node[rounded corners,rectangle,draw,fill=blue!20] (M) at (0,10) {\small $M$};
\node[rounded corners,rectangle,draw,fill=blue!20] (FLR) at (0,8) {\small $F_{LR}$};
\node[rounded corners,rectangle,draw,fill=blue!20] (FL) at (-3,6) {\small $F_L$};
\node[rounded corners,rectangle,draw,fill=blue!20] (GLR) at (0,6) {\small $G_{LR}$};
\node[rounded corners,rectangle,draw,fill=blue!20] (FR) at (3,6) {\small $F_R$};
\node[rounded corners,rectangle,draw,fill=blue!20] (GL) at (-3,4) {\small $G_L$};
\node[rounded corners,rectangle,draw,fill=blue!20] (F) at (0,4) {\small $F$};
\node[rounded corners,rectangle,draw,fill=blue!20] (GR) at (3,4) {\small $G_R$};
\node[rounded corners,rectangle,draw,fill=blue!20] (G) at (0,2) {\small $G$};
\node[rounded corners,rectangle,draw,fill=blue!20] (E) at (3,2) {\small $E$};
\node[rounded corners,rectangle,draw,fill=blue!20] (1) at (3,0) {\small $1$};
\end{scope}
\begin{scope}[shift={(10,0)}]
\node[rounded corners,rectangle,draw,fill=blue!20] (M) at (0,10) {};
\node[rounded corners,rectangle,draw,fill=blue!20] (FLR) at (0,8) {};
\node[rounded corners,rectangle,draw,fill=blue!20] (FL) at (-3,6) {};
\node[rounded corners,rectangle,draw,fill=blue!20] (GLR) at (0,6) {};
\node[rounded corners,rectangle,draw,fill=blue!20] (FR) at (3,6) {};
\node[rounded corners,rectangle,draw,fill=blue!20] (GL) at (-3,4) {};
\node[rounded corners,rectangle,draw,fill=blue!20] (F) at (0,4){};
\node[rounded corners,rectangle,draw,fill=blue!20] (GR) at (3,4) {};
\node[rounded corners,rectangle,draw,fill=blue!20] (G) at (0,2) {};
\draw (FLR)--(FL)--(GL)--(G)--(GR)--(FR)--(FLR);
\draw (GL)--(GLR)--(GR);
\draw[line width=1.5mm,white] (FL)--(F)--(FR);
\draw (FL)--(F)--(FR);
\draw (G)--(F);
\draw (M)--(FLR)--(GLR);
\node[rounded corners,rectangle,draw,fill=blue!20] (M) at (0,10) {\small $M$};
\node[rounded corners,rectangle,draw,fill=blue!20] (FLR) at (0,8) {\small $F_{LR}$};
\node[rounded corners,rectangle,draw,fill=blue!20] (FL) at (-3,6) {\small $F_L$};
\node[rounded corners,rectangle,draw,fill=blue!20] (GLR) at (0,6) {\small $G_{LR}$};
\node[rounded corners,rectangle,draw,fill=blue!20] (FR) at (3,6) {\small $F_R$};
\node[rounded corners,rectangle,draw,fill=blue!20] (GL) at (-3,4) {\small $G_L$};
\node[rounded corners,rectangle,draw,fill=blue!20] (F) at (0,4) {\small $F=E$};
\node[rounded corners,rectangle,draw,fill=blue!20] (GR) at (3,4) {\small $G_R$};
\node[rounded corners,rectangle,draw,fill=blue!20] (G) at (0,2) {\small $G=1$};
\node () at (0,11.2) {$(0,0,0,1)$};
\end{scope}
\begin{scope}[shift={(20,0)}]
\node[rounded corners,rectangle,draw,fill=blue!20] (M) at (0,10) {};
\node[rounded corners,rectangle,draw,fill=blue!20] (FLR) at (0,8) {};
\node[rounded corners,rectangle,draw,fill=blue!20] (FL) at (-3,6) {};
\node[rounded corners,rectangle,draw,fill=blue!20] (FR) at (3,6) {};
\node[rounded corners,rectangle,draw,fill=blue!20] (GL) at (-3,4) {};
\node[rounded corners,rectangle,draw,fill=blue!20] (F) at (0,4) {};
\node[rounded corners,rectangle,draw,fill=blue!20] (GR) at (3,4) {};
\node[rounded corners,rectangle,draw,fill=blue!20] (G) at (0,2) {};
\node[rounded corners,rectangle,draw,fill=blue!20] (E) at (3,2) {};
\node[rounded corners,rectangle,draw,fill=blue!20] (1) at (3,0) {};
\draw (FLR)--(FL)--(GL)--(G)--(GR)--(FR)--(FLR);
\draw[line width=1.5mm,white] (FL)--(F)--(FR);
\draw (FL)--(F)--(FR);
\draw (G)--(F) (E)--(1);
\draw (M)--(FLR);
\draw[line width=1.5mm,white] (F)--(E);
\draw (F)--(E) (G)--(1);
\node[rounded corners,rectangle,draw,fill=blue!20] (M) at (0,10) {\small $M$};
\node[rounded corners,rectangle,draw,fill=blue!20] (FLR) at (0,8) {\small $F_{LR}=G_{LR}$};
\node[rounded corners,rectangle,draw,fill=blue!20] (FL) at (-3,6) {\small $F_L$};
\node[rounded corners,rectangle,draw,fill=blue!20] (FR) at (3,6) {\small $F_R$};
\node[rounded corners,rectangle,draw,fill=blue!20] (GL) at (-3,4) {\small $G_L$};
\node[rounded corners,rectangle,draw,fill=blue!20] (F) at (0,4) {\small $F$};
\node[rounded corners,rectangle,draw,fill=blue!20] (GR) at (3,4) {\small $G_R$};
\node[rounded corners,rectangle,draw,fill=blue!20] (G) at (0,2) {\small $G$};
\node[rounded corners,rectangle,draw,fill=blue!20] (E) at (3,2) {\small $E$};
\node[rounded corners,rectangle,draw,fill=blue!20] (1) at (3,0) {\small $1$};
\node () at (0,11.2) {$(0,0,1,0)$};
\end{scope}
\begin{scope}[shift={(30,0)}]
\node[rounded corners,rectangle,draw,fill=blue!20] (M) at (0,10) {};
\node[rounded corners,rectangle,draw,fill=blue!20] (FLR) at (0,8) {};
\node[rounded corners,rectangle,draw,fill=blue!20] (FL) at (-3,6) {};
\node[rounded corners,rectangle,draw,fill=blue!20] (FR) at (3,6) {};
\node[rounded corners,rectangle,draw,fill=blue!20] (GL) at (-3,4) {};
\node[rounded corners,rectangle,draw,fill=blue!20] (F) at (0,4) {};
\node[rounded corners,rectangle,draw,fill=blue!20] (GR) at (3,4) {};
\node[rounded corners,rectangle,draw,fill=blue!20] (G) at (0,2) {};
\draw (FLR)--(FL)--(GL)--(G)--(GR)--(FR)--(FLR);
\draw[line width=1.5mm,white] (FL)--(F)--(FR);
\draw (FL)--(F)--(FR);
\draw (G)--(F);
\draw (M)--(FLR);
\node[rounded corners,rectangle,draw,fill=blue!20] (M) at (0,10) {\small $M$};
\node[rounded corners,rectangle,draw,fill=blue!20] (FLR) at (0,8) {\small $F_{LR}=G_{LR}$};
\node[rounded corners,rectangle,draw,fill=blue!20] (FL) at (-3,6) {\small $F_L$};
\node[rounded corners,rectangle,draw,fill=blue!20] (FR) at (3,6) {\small $F_R$};
\node[rounded corners,rectangle,draw,fill=blue!20] (GL) at (-3,4) {\small $G_L$};
\node[rounded corners,rectangle,draw,fill=blue!20] (F) at (0,4) {\small $F=E$};
\node[rounded corners,rectangle,draw,fill=blue!20] (GR) at (3,4) {\small $G_R$};
\node[rounded corners,rectangle,draw,fill=blue!20] (G) at (0,2) {\small $G=1$};
\node () at (0,11.2) {$(0,0,1,1)$};
\end{scope}
\begin{scope}[shift={(0,-12)}]
\node[rounded corners,rectangle,draw,fill=blue!20] (FLR) at (0,8) {};
\node[rounded corners,rectangle,draw,fill=blue!20] (FL) at (-3,6) {};
\node[rounded corners,rectangle,draw,fill=blue!20] (GLR) at (0,6) {};
\node[rounded corners,rectangle,draw,fill=blue!20] (FR) at (3,6) {};
\node[rounded corners,rectangle,draw,fill=blue!20] (GL) at (-3,4) {};
\node[rounded corners,rectangle,draw,fill=blue!20] (F) at (0,4) {};
\node[rounded corners,rectangle,draw,fill=blue!20] (GR) at (3,4) {};
\node[rounded corners,rectangle,draw,fill=blue!20] (G) at (0,2) {};
\node[rounded corners,rectangle,draw,fill=blue!20] (E) at (3,2) {};
\node[rounded corners,rectangle,draw,fill=blue!20] (1) at (3,0) {};
\draw (FLR)--(FL)--(GL)--(G)--(GR)--(FR)--(FLR);
\draw (GL)--(GLR)--(GR);
\draw[line width=1.5mm,white] (FL)--(F)--(FR);
\draw (FL)--(F)--(FR);
\draw (G)--(F) (E)--(1);
\draw (FLR)--(GLR);
\draw[line width=1.5mm,white] (F)--(E);
\draw (F)--(E) (G)--(1);
\node[rounded corners,rectangle,draw,fill=blue!20] (FLR) at (0,8) {\small $M=F_{LR}$};
\node[rounded corners,rectangle,draw,fill=blue!20] (FL) at (-3,6) {\small $F_L$};
\node[rounded corners,rectangle,draw,fill=blue!20] (GLR) at (0,6) {\small $G_{LR}$};
\node[rounded corners,rectangle,draw,fill=blue!20] (FR) at (3,6) {\small $F_R$};
\node[rounded corners,rectangle,draw,fill=blue!20] (GL) at (-3,4) {\small $G_L$};
\node[rounded corners,rectangle,draw,fill=blue!20] (F) at (0,4) {\small $F$};
\node[rounded corners,rectangle,draw,fill=blue!20] (GR) at (3,4) {\small $G_R$};
\node[rounded corners,rectangle,draw,fill=blue!20] (G) at (0,2) {\small $G$};
\node[rounded corners,rectangle,draw,fill=blue!20] (E) at (3,2) {\small $E$};
\node[rounded corners,rectangle,draw,fill=blue!20] (1) at (3,0) {\small $1$};
\node () at (0,9.2) {$(0,1,0,0)$};
\end{scope}
\begin{scope}[shift={(10,-12)}]
\node[rounded corners,rectangle,draw,fill=blue!20] (FLR) at (0,8) {};
\node[rounded corners,rectangle,draw,fill=blue!20] (FL) at (-3,6) {};
\node[rounded corners,rectangle,draw,fill=blue!20] (GLR) at (0,6) {};
\node[rounded corners,rectangle,draw,fill=blue!20] (FR) at (3,6) {};
\node[rounded corners,rectangle,draw,fill=blue!20] (GL) at (-3,4) {};
\node[rounded corners,rectangle,draw,fill=blue!20] (F) at (0,4) {};
\node[rounded corners,rectangle,draw,fill=blue!20] (GR) at (3,4) {};
\node[rounded corners,rectangle,draw,fill=blue!20] (G) at (0,2) {};
\draw (FLR)--(FL)--(GL)--(G)--(GR)--(FR)--(FLR);
\draw (GL)--(GLR)--(GR);
\draw[line width=1.5mm,white] (FL)--(F)--(FR);
\draw (FL)--(F)--(FR);
\draw (G)--(F);
\draw (FLR)--(GLR);
\node[rounded corners,rectangle,draw,fill=blue!20] (FLR) at (0,8) {\small $M=F_{LR}$};
\node[rounded corners,rectangle,draw,fill=blue!20] (FL) at (-3,6) {\small $F_L$};
\node[rounded corners,rectangle,draw,fill=blue!20] (GLR) at (0,6) {\small $G_{LR}$};
\node[rounded corners,rectangle,draw,fill=blue!20] (FR) at (3,6) {\small $F_R$};
\node[rounded corners,rectangle,draw,fill=blue!20] (GL) at (-3,4) {\small $G_L$};
\node[rounded corners,rectangle,draw,fill=blue!20] (F) at (0,4) {\small $F=E$};
\node[rounded corners,rectangle,draw,fill=blue!20] (GR) at (3,4) {\small $G_R$};
\node[rounded corners,rectangle,draw,fill=blue!20] (G) at (0,2) {\small $G=1$};
\node () at (0,9.2) {$(0,1,0,1)$};
\end{scope}
\begin{scope}[shift={(20,-12)}]
\node[rounded corners,rectangle,draw,fill=blue!20] (FLR) at (0,8) {};
\node[rounded corners,rectangle,draw,fill=blue!20] (FL) at (-3,6) {};
\node[rounded corners,rectangle,draw,fill=blue!20] (FR) at (3,6) {};
\node[rounded corners,rectangle,draw,fill=blue!20] (GL) at (-3,4) {};
\node[rounded corners,rectangle,draw,fill=blue!20] (F) at (0,4) {};
\node[rounded corners,rectangle,draw,fill=blue!20] (GR) at (3,4) {};
\node[rounded corners,rectangle,draw,fill=blue!20] (G) at (0,2) {};
\node[rounded corners,rectangle,draw,fill=blue!20] (E) at (3,2) {};
\node[rounded corners,rectangle,draw,fill=blue!20] (1) at (3,0) {};
\draw (FLR)--(FL)--(GL)--(G)--(GR)--(FR)--(FLR);
\draw[line width=1.5mm,white] (FL)--(F)--(FR);
\draw (FL)--(F)--(FR);
\draw (G)--(F) (E)--(1);
\draw[line width=1.5mm,white] (F)--(E);
\draw (F)--(E) (G)--(1);
\node[rounded corners,rectangle,draw,fill=blue!20] (FLR) at (0,8) {\small $M=F_{LR}=G_{LR}$};
\node[rounded corners,rectangle,draw,fill=blue!20] (FL) at (-3,6) {\small $F_L$};
\node[rounded corners,rectangle,draw,fill=blue!20] (FR) at (3,6) {\small $F_R$};
\node[rounded corners,rectangle,draw,fill=blue!20] (GL) at (-3,4) {\small $G_L$};
\node[rounded corners,rectangle,draw,fill=blue!20] (F) at (0,4) {\small $F$};
\node[rounded corners,rectangle,draw,fill=blue!20] (GR) at (3,4) {\small $G_R$};
\node[rounded corners,rectangle,draw,fill=blue!20] (G) at (0,2) {\small $G$};
\node[rounded corners,rectangle,draw,fill=blue!20] (E) at (3,2) {\small $E$};
\node[rounded corners,rectangle,draw,fill=blue!20] (1) at (3,0) {\small $1$};
\node () at (0,9.2) {$(0,1,1,0)$};
\end{scope}
\begin{scope}[shift={(30,-12)}]
\node[rounded corners,rectangle,draw,fill=blue!20] (FLR) at (0,8) {};
\node[rounded corners,rectangle,draw,fill=blue!20] (FL) at (-3,6) {};
\node[rounded corners,rectangle,draw,fill=blue!20] (FR) at (3,6) {};
\node[rounded corners,rectangle,draw,fill=blue!20] (GL) at (-3,4) {};
\node[rounded corners,rectangle,draw,fill=blue!20] (F) at (0,4) {};
\node[rounded corners,rectangle,draw,fill=blue!20] (GR) at (3,4) {};
\node[rounded corners,rectangle,draw,fill=blue!20] (G) at (0,2) {};
\draw (FLR)--(FL)--(GL)--(G)--(GR)--(FR)--(FLR);
\draw[line width=1.5mm,white] (FL)--(F)--(FR);
\draw (FL)--(F)--(FR);
\draw (G)--(F);
\node[rounded corners,rectangle,draw,fill=blue!20] (FLR) at (0,8) {\small $M=F_{LR}=G_{LR}$};
\node[rounded corners,rectangle,draw,fill=blue!20] (FL) at (-3,6) {\small $F_L$};
\node[rounded corners,rectangle,draw,fill=blue!20] (FR) at (3,6) {\small $F_R$};
\node[rounded corners,rectangle,draw,fill=blue!20] (GL) at (-3,4) {\small $G_L$};
\node[rounded corners,rectangle,draw,fill=blue!20] (F) at (0,4) {\small $F=E$};
\node[rounded corners,rectangle,draw,fill=blue!20] (GR) at (3,4) {\small $G_R$};
\node[rounded corners,rectangle,draw,fill=blue!20] (G) at (0,2) {\small $G=1$};
\node () at (0,9.2) {$(0,1,1,1)$};
\end{scope}
\end{tikzpicture}
}
\caption{The lattice $\L(M)$ when $M$ has an unstable identity, according to the type $\T(M)=(0,i,j,k)$.  In each case, the nodes represent distinct submonoids of $M$.}
\label{f:L2}
\end{center}
\end{figure}

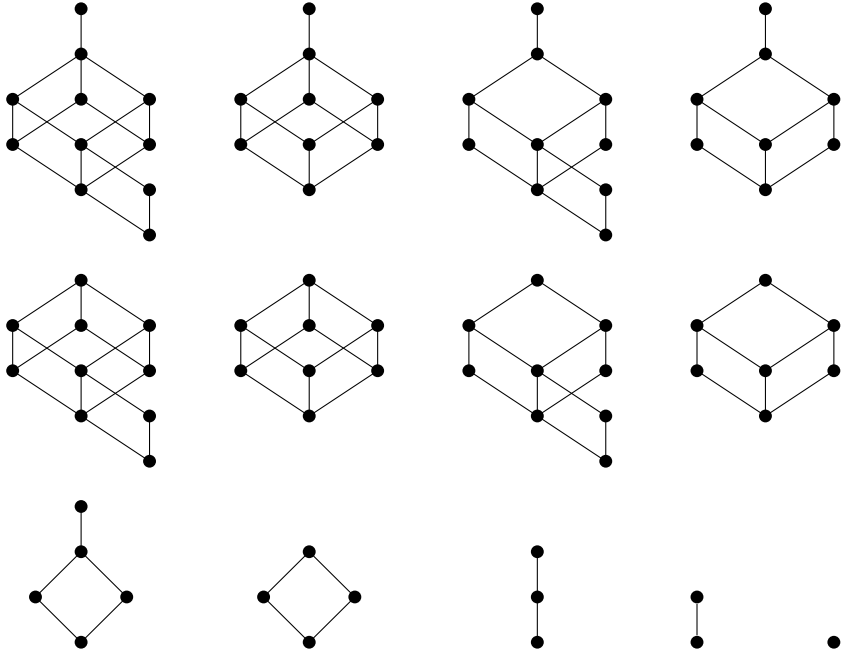
\begin{figure}[ht]
\begin{center}
\begin{tikzpicture}[scale=0.3]
\begin{scope}[shift={(0,0)}]
\node[circle,fill=black, inner sep = 0.0cm] (FLR) at (0,8) {};
\node[circle,fill=black, inner sep = 0.0cm] (FL) at (-3,6) {};
\node[circle,fill=black, inner sep = 0.0cm] (FR) at (3,6) {};
\node[circle,fill=black, inner sep = 0.0cm] (F) at (0,4) {};
\node[circle,fill=black, inner sep = 0.0cm] (GL) at (-3,4) {};
\node[circle,fill=black, inner sep = 0.0cm] (GR) at (3,4) {};
\node[circle,fill=black, inner sep = 0.0cm] (G) at (0,2) {};
\draw (G)--(GR)--(FR)--(FLR)--(FL)--(GL)--(G) (FL)--(F)--(FR) (G)--(F);
\node[circle,fill=black, inner sep = 0.0cm] (E) at (3,2) {}; \node[circle,fill=black, inner sep = 0.0cm] (O) at (3,0) {}; \draw (F)--(E)--(O)--(G);
\node[circle,fill=black, inner sep = 0.0cm] (GLR) at (0,6) {}; \draw (GL)--(GLR)--(GR) (FLR)--(GLR);
\node[circle,fill=black, inner sep = 0.0cm] (M) at (0,10) {}; \draw(M)--(FLR);
\node[circle,fill=black, inner sep = 0.06cm] (FLR) at (0,8) {};
\node[circle,fill=black, inner sep = 0.06cm] (FL) at (-3,6) {};
\node[circle,fill=black, inner sep = 0.06cm] (FR) at (3,6) {};
\node[circle,fill=black, inner sep = 0.06cm] (F) at (0,4) {};
\node[circle,fill=black, inner sep = 0.06cm] (GL) at (-3,4) {};
\node[circle,fill=black, inner sep = 0.06cm] (GR) at (3,4) {};
\node[circle,fill=black, inner sep = 0.06cm] (G) at (0,2) {};
\node[circle,fill=black, inner sep = 0.06cm] (E) at (3,2) {}; \node[circle,fill=black, inner sep = 0.06cm] (O) at (3,0) {}; 
\node[circle,fill=black, inner sep = 0.06cm] (GLR) at (0,6) {}; 
\node[circle,fill=black, inner sep = 0.06cm] (M) at (0,10) {}; 
\end{scope}
\begin{scope}[shift={(10,0)}]
\node[circle,fill=black, inner sep = 0.0cm] (FLR) at (0,8) {};
\node[circle,fill=black, inner sep = 0.0cm] (FL) at (-3,6) {};
\node[circle,fill=black, inner sep = 0.0cm] (FR) at (3,6) {};
\node[circle,fill=black, inner sep = 0.0cm] (F) at (0,4) {};
\node[circle,fill=black, inner sep = 0.0cm] (GL) at (-3,4) {};
\node[circle,fill=black, inner sep = 0.0cm] (GR) at (3,4) {};
\node[circle,fill=black, inner sep = 0.0cm] (G) at (0,2) {};
\draw (G)--(GR)--(FR)--(FLR)--(FL)--(GL)--(G) (FL)--(F)--(FR) (G)--(F);
\node[circle,fill=black, inner sep = 0.0cm] (GLR) at (0,6) {}; \draw (GL)--(GLR)--(GR) (FLR)--(GLR);
\node[circle,fill=black, inner sep = 0.0cm] (M) at (0,10) {}; \draw(M)--(FLR);
\node[circle,fill=black, inner sep = 0.06cm] (FLR) at (0,8) {};
\node[circle,fill=black, inner sep = 0.06cm] (FL) at (-3,6) {};
\node[circle,fill=black, inner sep = 0.06cm] (FR) at (3,6) {};
\node[circle,fill=black, inner sep = 0.06cm] (F) at (0,4) {};
\node[circle,fill=black, inner sep = 0.06cm] (GL) at (-3,4) {};
\node[circle,fill=black, inner sep = 0.06cm] (GR) at (3,4) {};
\node[circle,fill=black, inner sep = 0.06cm] (G) at (0,2) {};
\node[circle,fill=black, inner sep = 0.06cm] (GLR) at (0,6) {}; 
\node[circle,fill=black, inner sep = 0.06cm] (M) at (0,10) {}; 
\end{scope}
\begin{scope}[shift={(20,0)}]
\node[circle,fill=black, inner sep = 0.0cm] (FLR) at (0,8) {};
\node[circle,fill=black, inner sep = 0.0cm] (FL) at (-3,6) {};
\node[circle,fill=black, inner sep = 0.0cm] (FR) at (3,6) {};
\node[circle,fill=black, inner sep = 0.0cm] (F) at (0,4) {};
\node[circle,fill=black, inner sep = 0.0cm] (GL) at (-3,4) {};
\node[circle,fill=black, inner sep = 0.0cm] (GR) at (3,4) {};
\node[circle,fill=black, inner sep = 0.0cm] (G) at (0,2) {};
\draw (G)--(GR)--(FR)--(FLR)--(FL)--(GL)--(G) (FL)--(F)--(FR) (G)--(F);
\node[circle,fill=black, inner sep = 0.0cm] (E) at (3,2) {}; \node[circle,fill=black, inner sep = 0.0cm] (O) at (3,0) {}; \draw (F)--(E)--(O)--(G);
\node[circle,fill=black, inner sep = 0.0cm] (M) at (0,10) {}; \draw(M)--(FLR);
\node[circle,fill=black, inner sep = 0.06cm] (FLR) at (0,8) {};
\node[circle,fill=black, inner sep = 0.06cm] (FL) at (-3,6) {};
\node[circle,fill=black, inner sep = 0.06cm] (FR) at (3,6) {};
\node[circle,fill=black, inner sep = 0.06cm] (F) at (0,4) {};
\node[circle,fill=black, inner sep = 0.06cm] (GL) at (-3,4) {};
\node[circle,fill=black, inner sep = 0.06cm] (GR) at (3,4) {};
\node[circle,fill=black, inner sep = 0.06cm] (G) at (0,2) {};
\node[circle,fill=black, inner sep = 0.06cm] (E) at (3,2) {}; \node[circle,fill=black, inner sep = 0.06cm] (O) at (3,0) {}; 
\node[circle,fill=black, inner sep = 0.06cm] (M) at (0,10) {}; 
\end{scope}
\begin{scope}[shift={(30,0)}]
\node[circle,fill=black, inner sep = 0.0cm] (FLR) at (0,8) {};
\node[circle,fill=black, inner sep = 0.0cm] (FL) at (-3,6) {};
\node[circle,fill=black, inner sep = 0.0cm] (FR) at (3,6) {};
\node[circle,fill=black, inner sep = 0.0cm] (F) at (0,4) {};
\node[circle,fill=black, inner sep = 0.0cm] (GL) at (-3,4) {};
\node[circle,fill=black, inner sep = 0.0cm] (GR) at (3,4) {};
\node[circle,fill=black, inner sep = 0.0cm] (G) at (0,2) {};
\draw (G)--(GR)--(FR)--(FLR)--(FL)--(GL)--(G) (FL)--(F)--(FR) (G)--(F);
\node[circle,fill=black, inner sep = 0.0cm] (M) at (0,10) {}; \draw(M)--(FLR);
\node[circle,fill=black, inner sep = 0.06cm] (FLR) at (0,8) {};
\node[circle,fill=black, inner sep = 0.06cm] (FL) at (-3,6) {};
\node[circle,fill=black, inner sep = 0.06cm] (FR) at (3,6) {};
\node[circle,fill=black, inner sep = 0.06cm] (F) at (0,4) {};
\node[circle,fill=black, inner sep = 0.06cm] (GL) at (-3,4) {};
\node[circle,fill=black, inner sep = 0.06cm] (GR) at (3,4) {};
\node[circle,fill=black, inner sep = 0.06cm] (G) at (0,2) {};
\node[circle,fill=black, inner sep = 0.06cm] (M) at (0,10) {}; 
\end{scope}
\begin{scope}[shift={(0,-10)}]
\node[circle,fill=black, inner sep = 0.0cm] (FLR) at (0,8) {};
\node[circle,fill=black, inner sep = 0.0cm] (FL) at (-3,6) {};
\node[circle,fill=black, inner sep = 0.0cm] (FR) at (3,6) {};
\node[circle,fill=black, inner sep = 0.0cm] (F) at (0,4) {};
\node[circle,fill=black, inner sep = 0.0cm] (GL) at (-3,4) {};
\node[circle,fill=black, inner sep = 0.0cm] (GR) at (3,4) {};
\node[circle,fill=black, inner sep = 0.0cm] (G) at (0,2) {};
\draw (G)--(GR)--(FR)--(FLR)--(FL)--(GL)--(G) (FL)--(F)--(FR) (G)--(F);
\node[circle,fill=black, inner sep = 0.0cm] (E) at (3,2) {}; \node[circle,fill=black, inner sep = 0.0cm] (O) at (3,0) {}; \draw (F)--(E)--(O)--(G);
\node[circle,fill=black, inner sep = 0.0cm] (GLR) at (0,6) {}; \draw (GL)--(GLR)--(GR) (FLR)--(GLR);
\node[circle,fill=black, inner sep = 0.06cm] (FLR) at (0,8) {};
\node[circle,fill=black, inner sep = 0.06cm] (FL) at (-3,6) {};
\node[circle,fill=black, inner sep = 0.06cm] (FR) at (3,6) {};
\node[circle,fill=black, inner sep = 0.06cm] (F) at (0,4) {};
\node[circle,fill=black, inner sep = 0.06cm] (GL) at (-3,4) {};
\node[circle,fill=black, inner sep = 0.06cm] (GR) at (3,4) {};
\node[circle,fill=black, inner sep = 0.06cm] (G) at (0,2) {};
\node[circle,fill=black, inner sep = 0.06cm] (E) at (3,2) {}; \node[circle,fill=black, inner sep = 0.06cm] (O) at (3,0) {}; 
\node[circle,fill=black, inner sep = 0.06cm] (GLR) at (0,6) {}; 
\end{scope}
\begin{scope}[shift={(10,-10)}]
\node[circle,fill=black, inner sep = 0.0cm] (FLR) at (0,8) {};
\node[circle,fill=black, inner sep = 0.0cm] (FL) at (-3,6) {};
\node[circle,fill=black, inner sep = 0.0cm] (FR) at (3,6) {};
\node[circle,fill=black, inner sep = 0.0cm] (F) at (0,4) {};
\node[circle,fill=black, inner sep = 0.0cm] (GL) at (-3,4) {};
\node[circle,fill=black, inner sep = 0.0cm] (GR) at (3,4) {};
\node[circle,fill=black, inner sep = 0.0cm] (G) at (0,2) {};
\draw (G)--(GR)--(FR)--(FLR)--(FL)--(GL)--(G) (FL)--(F)--(FR) (G)--(F);
\node[circle,fill=black, inner sep = 0.0cm] (GLR) at (0,6) {}; \draw (GL)--(GLR)--(GR) (FLR)--(GLR);
\node[circle,fill=black, inner sep = 0.06cm] (FLR) at (0,8) {};
\node[circle,fill=black, inner sep = 0.06cm] (FL) at (-3,6) {};
\node[circle,fill=black, inner sep = 0.06cm] (FR) at (3,6) {};
\node[circle,fill=black, inner sep = 0.06cm] (F) at (0,4) {};
\node[circle,fill=black, inner sep = 0.06cm] (GL) at (-3,4) {};
\node[circle,fill=black, inner sep = 0.06cm] (GR) at (3,4) {};
\node[circle,fill=black, inner sep = 0.06cm] (G) at (0,2) {};
\node[circle,fill=black, inner sep = 0.06cm] (GLR) at (0,6) {}; 
\end{scope}
\begin{scope}[shift={(20,-10)}]
\node[circle,fill=black, inner sep = 0.0cm] (FLR) at (0,8) {};
\node[circle,fill=black, inner sep = 0.0cm] (FL) at (-3,6) {};
\node[circle,fill=black, inner sep = 0.0cm] (FR) at (3,6) {};
\node[circle,fill=black, inner sep = 0.0cm] (F) at (0,4) {};
\node[circle,fill=black, inner sep = 0.0cm] (GL) at (-3,4) {};
\node[circle,fill=black, inner sep = 0.0cm] (GR) at (3,4) {};
\node[circle,fill=black, inner sep = 0.0cm] (G) at (0,2) {};
\draw (G)--(GR)--(FR)--(FLR)--(FL)--(GL)--(G) (FL)--(F)--(FR) (G)--(F);
\node[circle,fill=black, inner sep = 0.0cm] (E) at (3,2) {}; \node[circle,fill=black, inner sep = 0.0cm] (O) at (3,0) {}; \draw (F)--(E)--(O)--(G);
\node[circle,fill=black, inner sep = 0.06cm] (FLR) at (0,8) {};
\node[circle,fill=black, inner sep = 0.06cm] (FL) at (-3,6) {};
\node[circle,fill=black, inner sep = 0.06cm] (FR) at (3,6) {};
\node[circle,fill=black, inner sep = 0.06cm] (F) at (0,4) {};
\node[circle,fill=black, inner sep = 0.06cm] (GL) at (-3,4) {};
\node[circle,fill=black, inner sep = 0.06cm] (GR) at (3,4) {};
\node[circle,fill=black, inner sep = 0.06cm] (G) at (0,2) {};
\node[circle,fill=black, inner sep = 0.06cm] (E) at (3,2) {}; \node[circle,fill=black, inner sep = 0.06cm] (O) at (3,0) {}; 
\end{scope}
\begin{scope}[shift={(30,-10)}]
\node[circle,fill=black, inner sep = 0.0cm] (FLR) at (0,8) {};
\node[circle,fill=black, inner sep = 0.0cm] (FL) at (-3,6) {};
\node[circle,fill=black, inner sep = 0.0cm] (FR) at (3,6) {};
\node[circle,fill=black, inner sep = 0.0cm] (F) at (0,4) {};
\node[circle,fill=black, inner sep = 0.0cm] (GL) at (-3,4) {};
\node[circle,fill=black, inner sep = 0.0cm] (GR) at (3,4) {};
\node[circle,fill=black, inner sep = 0.0cm] (G) at (0,2) {};
\draw (G)--(GR)--(FR)--(FLR)--(FL)--(GL)--(G) (FL)--(F)--(FR) (G)--(F);
\node[circle,fill=black, inner sep = 0.06cm] (FLR) at (0,8) {};
\node[circle,fill=black, inner sep = 0.06cm] (FL) at (-3,6) {};
\node[circle,fill=black, inner sep = 0.06cm] (FR) at (3,6) {};
\node[circle,fill=black, inner sep = 0.06cm] (F) at (0,4) {};
\node[circle,fill=black, inner sep = 0.06cm] (GL) at (-3,4) {};
\node[circle,fill=black, inner sep = 0.06cm] (GR) at (3,4) {};
\node[circle,fill=black, inner sep = 0.06cm] (G) at (0,2) {};
\end{scope}
\begin{scope}[shift={(0,-22)}]
\node[circle,fill=black, inner sep = 0.0cm] (M) at (0,10) {};
\node[circle,fill=black, inner sep = 0.0cm] (F) at (0,8) {};
\node[circle,fill=black, inner sep = 0.0cm] (G) at (-2,6) {};
\node[circle,fill=black, inner sep = 0.0cm] (E) at (2,6) {};
\node[circle,fill=black, inner sep = 0.0cm] (1) at (0,4) {};
\draw (M)--(F);
\draw (F)--(E) (G)--(1);
\draw (F)--(G) (E)--(1);
\node[circle,fill=black, inner sep = 0.06cm] (M) at (0,10) {};
\node[circle,fill=black, inner sep = 0.06cm] (F) at (0,8) {};
\node[circle,fill=black, inner sep = 0.06cm] (G) at (-2,6) {};
\node[circle,fill=black, inner sep = 0.06cm] (E) at (2,6) {};
\node[circle,fill=black, inner sep = 0.06cm] (1) at (0,4) {};
\end{scope}
\begin{scope}[shift={(10,-22)}]
\node[circle,fill=black, inner sep = 0.0cm] (F) at (0,8) {};
\node[circle,fill=black, inner sep = 0.0cm] (G) at (-2,6) {};
\node[circle,fill=black, inner sep = 0.0cm] (E) at (2,6) {};
\node[circle,fill=black, inner sep = 0.0cm] (1) at (0,4) {};
\draw (F)--(E) (G)--(1);
\draw (F)--(G) (E)--(1);
\node[circle,fill=black, inner sep = 0.06cm] (F) at (0,8) {};
\node[circle,fill=black, inner sep = 0.06cm] (G) at (-2,6) {};
\node[circle,fill=black, inner sep = 0.06cm] (E) at (2,6) {};
\node[circle,fill=black, inner sep = 0.06cm] (1) at (0,4) {};
\end{scope}
\begin{scope}[shift={(20,-22)}]
\node[circle,fill=black, inner sep = 0.0cm] (M) at (0,8) {};
\node[circle,fill=black, inner sep = 0.0cm] (F) at (0,6) {};
\node[circle,fill=black, inner sep = 0.0cm] (1) at (0,4) {};
\draw (M)--(F);
\draw (F)--(1);
\node[circle,fill=black, inner sep = 0.06cm] (M) at (0,8) {};
\node[circle,fill=black, inner sep = 0.06cm] (F) at (0,6) {};
\node[circle,fill=black, inner sep = 0.06cm] (1) at (0,4) {};
\end{scope}
\begin{scope}[shift={(27,-22)}]
\node[circle,fill=black, inner sep = 0.06cm] (F) at (0,6) {};
\node[circle,fill=black, inner sep = 0.06cm] (1) at (0,4) {};
\draw (F)--(1);
\node[circle,fill=black, inner sep = 0.0cm] (F) at (0,6) {};
\node[circle,fill=black, inner sep = 0.0cm] (1) at (0,4) {};
\end{scope}
\begin{scope}[shift={(33,-22)}]
\node[circle,fill=black, inner sep = 0.06cm] (1) at (0,4) {};
\end{scope}
\end{tikzpicture}
\caption{The possible lattices $\L(M)$ for a monoid $M$, up to lattice isomorphism.}
\label{f:L3}
\end{center}
\end{figure}

In this section we classify the lattices $\L(M)$, for monoids $M$.  To do so, we first define the \emph{type} of a monoid, as a certain binary tuple of length $4$; we show in Proposition \ref{p:types} that all sixteen such tuples occur as the type of a monoid.  We then show that the type of $M$ completely determines the lattice $\L(M)$, considering separate cases in which the identity of $M$ is stable or unstable.  The classification is stated in Theorem \ref{t:L}; cf.~Figures \ref{f:lattice_small}--\ref{f:L3}.
Throughout this section, unless otherwise specified, $M$ denotes an arbitrary monoid, and we continue to use the abbreviations \eqref{e:G} and \eqref{e:F}.

Consider the following questions concerning a monoid $M$:
\benT
\item \label{T1} Does $G=G_L$ hold?
\item \label{T2} Does $F_{LR}=M$ hold?
\item \label{T3} Does $F_{LR}=G_{LR}$ hold?
\item \label{T4} Does $G=\{1\}$ hold?
\een
We denote the Yes (=1) or No~(=0) answers to these questions by $\Tone(M)$, $\Ttwo(M)$, $\Tthree(M)$ and $\Tfour(M)$, respectively.  We also define the binary quadruple
\[
\T(M)=(\Tone(M),\Ttwo(M),\Tthree(M),\Tfour(M)),
\]
and call this the \emph{type of $M$}.  There are sixteen quadruples over $\{0,1\}$, and Proposition \ref{p:types} below shows that each such quadruple is the type of some monoid.

By Lemma \ref{l:DP}, if $\X,\Y\in\Fu$ then for any monoids $M$ and $N$, we have
\[
\X(M\times N)=\Y(M\times N) \iff \X(M)=\Y(M)\text{ and } \X(N)=\Y(N).
\]
It follows that the integers $\Ti(M)$ are multiplicative, in the sense that for monoids $M$ and $N$, we have $\Ti(M\times N)=\Ti(M)\times\Ti(N)$; here the first $\times$ is monoid direct product, and the second is ordinary integer multiplication in $\{0,1\}$.  It follows that types are multiplicative as well:
\begin{equation}\label{e:T}
\T(M\times N) = \T(M)\times\T(N) \qquad\text{for monoids $M$ and $N$.}
\end{equation}
In the second expression, we mean the coordinate-wise product of tuples.

\begin{prop}\label{p:types}
For any $i,j,k,l\in\{0,1\}$, there exists a monoid $M$ with type $\T(M)=(i,j,k,l)$.
\end{prop}

\pf
Consulting Table \ref{t:GNB_small}, we see that
\bit
\item $\T(G)=(1,1,1,0)$ for a nontrivial group $G$,
\item $\T(E)=(1,1,0,1)$ for a nontrivial idempotent-generated monoid $E$,
\item $\T(\Pos)=(1,0,1,1)$ for the multiplicative monoid of positive integers $\Pos$,
\item $\T(B)=(0,1,1,1)$ for the bicyclic monoid $B$.
\eit
Thus, in light of \eqref{e:T}, we can obtain a monoid with any type by taking a suitable direct product of some (possibly empty) collection of $G$, $E$, $\Pos$, $B$, as above.  
\epf

The rest of Section \ref{s:class} is devoted to showing that the type of the monoid $M$ completely determines the lattice $\L(M)$.

We first consider the case in which the identity of $M$ is stable.  By Lemma \ref{l:J}, this is equivalent to having $G=G_L$: i.e., to having $\Tone(M)=1$.
In this case, the conditions in Lemma~\ref{l:J} all hold, but the conditions in Lemma \ref{l:GLR} do not (cf.~Remark \ref{r:dichotomy}).  
In particular, we have $G=G_L=G_R=G_{LR}$ and $F=F_L=F_R=F_{LR}$.  Thus, the lattice
\[
\L(M) = \big\{\{1\},E,G,F,M\big\}
\]
simplifies substantially, and has the generic shape pictured in Figure~\ref{f:lattice_small}.  In this diagram and others to follow, the trivial submonoid $\{1\}$ is abbreviated to $1$.

In general, some of the submonoids pictured in Figure \ref{f:lattice_small} could be equal, but by Lemma \ref{l:EG} (and the fact that $F=F_{LR}$) we have
\[
G=\{1\}\iff F=E \AND E=\{1\}\iff F=G.
\]
Also note that since $F=F_{LR}$ and $G=G_{LR}$, questions \ref{T2} and \ref{T3} are equivalent (in the case of $M$ having a stable identity) to:
\benTd\addtocounter{enumi}{1}
\item \label{T2'} Does $F=M$ hold?
\item \label{T3'} Does $F=G$ (equivalently, $E=\{1\}$) hold?
\een
Figure \ref{f:L1} shows the lattice $\L(M)$ for monoids of type $(1,i,j,k)$.  The values of $i=\Ttwod(M)$, $j=\Tthreed(M)$ and $k=\Tfour(M)$ determine which edges (if any) in Figure \ref{f:lattice_small} are contracted.

We now consider the case in which the identity of $M$ is unstable, which is equivalent to having $G\not=G_L$: i.e., to having $\Tone(M)=0$.
In this case, the conditions in Lemma~\ref{l:GLR} all hold, but the conditions in Lemma~\ref{l:J} do not.  
In particular, $G$, $G_L$, $G_R$ and $G_{LR}$ are four distinct submonoids; so too are $F$, $F_L$, $F_R$ and~$F_{LR}$.  Moreover, $G_{LR}$ (and hence $M$) contains infinitely many idempotents (cf.~Lemma \ref{l:GLR}), so certainly $E\not=1$; it follows from Lemma \ref{l:EG}\ref{EG1} that $\{G,G_L,G_R\} \cap \{F,F_L,F_R,F_{LR}\} \not=\emptyset$.  All of the above shows that the following seven submonoids of $M$ are distinct:
\begin{equation}\label{e:7}
\text{$G$, $G_L$, $G_R$, $F$, $F_L$, $F_R$ and $F_{LR}$.}
\end{equation}
These submonoids are shaded red in Figure \ref{f:lattice_red}, which gives the generic shape of $\L(M)$ in the unstable case.  Again we note that $G=\{1\}\iff F=E$; cf.~Lemma \ref{l:EG}\ref{EG2}.

Figure \ref{f:L2} shows the shapes the lattice $\L(M)$ takes for monoids of each type $(0,i,j,k)$, and again the values of $i,j,k$ determine which thin edges (if any) in Figure \ref{f:lattice_red} are contracted.

The results of this section may be summarised as follows:

\begin{thm}\label{t:L}
\ben
\item If a monoid $M$ has a stable identity, then the lattice $\L(M)$ is as shown in Figure \ref{f:L1}, according to its type $\T(M)=(1,i,j,k)$.
\item If a monoid $M$ has an unstable identity, then the lattice $\L(M)$ is as shown in Figure \ref{f:L2}, according to its type $\T(M)=(0,i,j,k)$.
\item Each of the lattices pictured in Figures \ref{f:L1} and \ref{f:L2} arises as $\L(M)$ for some monoid $M$.  
\item Up to isomorphism, the lattice $\L(M)$ associated to a monoid $M$ has one of the forms shown in Figure~\ref{f:L3}. \qed
\een
\end{thm}

\section{A semigroup of functors}\label{s:F}

In this section we study the semigroup of functors $\M\to\M$ generated (via composition) by the functors considered so far:
\[
\Fu = \{\O,\E,\G,\GL,\GR,\GLR,\F,\FL,\FR,\FLR,\I\}.
\]
We begin in Section \ref{ss:comp} by calculating compositions of the functors from $\Fu$, and observe that four such compositions do not seem to belong to $\Fu$.  In Section \ref{ss:more} we define a suitably enlarged set $\Fup$ of functors, and associate an enhanced lattice $\Lp(M)$ to each monoid $M$.  In Section \ref{ss:more2} we show that $\Fup$ is a semigroup, indeed a monoid; we calculate its size in Section \ref{ss:size}, and describe its algebraic structure in Section \ref{ss:structure}.  In Section~\ref{ss:LF} we calculate the lattice $\L(\Fup)$.

Throughout this section, unless otherwise specified, $M$ denotes an arbitrary monoid, and we continue to use the abbreviations \eqref{e:G} and \eqref{e:F}.

\subsection{Compositions}\label{ss:comp}

Since each functor from $\Fu$ maps $\M\to\M$, these functors may be composed.  For example, we may consider the functor $\E\circ\G:\M\to\M$.  Since groups have only one idempotent, we have $\E(\G(M))=\{1\}=\O(M)$ for any monoid~$M$, and this means that $\E\circ\G=\O$.  On the other hand, we have $\E\circ\E=\E$.
We also clearly have
\[
\O\circ\X=\X\circ\O=\O \AND \I\circ\X=\X\circ\I=\X \qquad\text{for any $\X\in\Fu$.}
\]
Various results from \cite[Section 2]{JE_IBM} may be interpreted as further such compositional equations.  
For example, \cite[Lemmas~2.1 and~2.9]{JE_IBM} say that if $\X$ is one of $\G$, $\GL$ or $\GR$, then
\[
\E\circ\X=\X\circ\E=\O
\AND
\G\circ\X = 
\GL\circ\X = 
\GR\circ\X = 
\F\circ\X = 
\FL\circ\X = 
\FR\circ\X = 
\G.
\]
Similarly, \cite[Lemma 2.8]{JE_IBM} says that if $\X$ is any of $\F$, $\FL$ or $\FR$, then
\[
\E\circ\X=\E \COMMA
\G\circ\X=\GL\circ\X=\GR\circ\X=\G \COMMA
\F\circ\X=\FL\circ\X=\FR\circ\X=\F.
\]
If $\hs$ represents any subscript other than $LR$, then since $\G_\hs\circ\E=\O$ (noted above), we have
\[
\F_\hs\circ\E(M) = \F_\hs(E) = \E(E)\vee\G_\hs(E) = \E(M)\vee\{1\} = \E(M),
\]
so that $\F_\hs\circ\E=\E$.  
The above composition rules are recorded as the black entries in Table \ref{t:op1}.  

\begin{table}[t]
\begin{center}
\begin{tabular}{l|lllllllllll}
\hspace{.5mm}$\circ$ & $\O$ & $\E$ & $\G$  & $\GL$  & $\GR$  & $\GLR$ & $\F$  & $\FL$  & $\FR$  & $\FLR$  & $\I$ \\
\hline
$\O$ & $\O$ & $\O$ & $\O$ & $\O$ & $\O$ & $\O$ & $\O$ & $\O$ & $\O$ & $\O$ & $\O$ \\
$\E$ & $\O$ & $\E$ & $\O$ & $\O$ & $\O$ &  & $\E$ & $\E$ & $\E$ & \textcolor{cyan}{$\E$} & $\E$ \\
$\G$ & $\O$ & $\O$ & $\G$ & $\G$ & $\G$ & \textcolor{cyan}{ $\G$} & $\G$ & $\G$ & $\G$ & \textcolor{cyan}{ $\G$} & $\G$ \\
$\GL$ & $\O$ & $\O$ & $\G$ & $\G$ & $\G$ & \textcolor{cyan}{ $\GL$} & $\G$ & $\G$ & $\G$ & \textcolor{cyan}{ $\GL$} & $\GL$ \\
$\GR$ & $\O$ & $\O$ & $\G$ & $\G$ & $\G$ & \textcolor{cyan}{ $\GR$} & $\G$ & $\G$ & $\G$ & \textcolor{cyan}{ $\GR$} & $\GR$ \\
$\GLR$ & $\O$ & \textcolor{cyan}{ $\O$} & \textcolor{cyan}{ $\G$} & \textcolor{cyan}{ $\G$} & \textcolor{cyan}{ $\G$} & \textcolor{cyan}{ $\GLR$} & \textcolor{cyan}{ $\G$} & \textcolor{cyan}{ $\G$} & \textcolor{cyan}{ $\G$} & \textcolor{cyan}{ $\GLR$} & $\GLR$ \\
$\F$ & $\O$ & $\E$ & $\G$ & $\G$ & $\G$ &  & $\F$ & $\F$ & $\F$ & \textcolor{cyan}{ $\F$} & $\F$ \\
$\FL$ & $\O$ & $\E$ & $\G$ & $\G$ & $\G$ &  & $\F$ & $\F$ & $\F$ & \textcolor{cyan}{ $\FL$} & $\FL$ \\
$\FR$ & $\O$ & $\E$ & $\G$ & $\G$ & $\G$ &  & $\F$ & $\F$ & $\F$ & \textcolor{cyan}{ $\FR$} & $\FR$ \\
$\FLR$ & $\O$ & \textcolor{cyan}{ $\E$} & \textcolor{cyan}{ $\G$} & \textcolor{cyan}{ $\G$} & \textcolor{cyan}{ $\G$} & \textcolor{cyan}{ $\GLR$} & \textcolor{cyan}{ $\F$} & \textcolor{cyan}{ $\F$} & \textcolor{cyan}{ $\F$} & \textcolor{cyan}{ $\FLR$} & $\FLR$ \\
$\I$ &  $\O$ & $\E$ & $\G$  & $\G_L$  & $\G_R$  & $\G_{LR}$ & $\F$  & $\FL$  & $\FR$  & $\FLR$  & $\I$ \\
\end{tabular}
\caption{Composition of the functors from $\Fu$.}
\label{t:op1}
\end{center}
\end{table}

Table \ref{t:op1} contains a number of other entries in blue (and some missing entries, which we will discuss in Section \ref{ss:more}).  The blue entries follow from Lemmas \ref{l:G} and \ref{l:F} below.  The proofs of these will use the following simple fact.

\begin{lemma}\label{l:X}
If $N$ is a submonoid of $M$, and if $\GL(M),\GR(M)\sub N$, then 
\[
\X(N)=\X(M) \qquad\text{for $\X=\G,\GL,\GR,\GLR$.}
\]
\end{lemma}

\pf
We first prove the claim for $\X=\GL$.  Since $N\sub M$, we clearly have $\GL(N)\sub\GL(M)$.  Conversely, suppose $x\in \GL(M)$.  So $x\in N$ by assumption.  We also have $1=ax$ for some $a\in M$.  But this implies that $a\in\GR(M)\sub N$, so in fact $x\in\GL(N)$ as required.  

The claim for $\X=\GR$ is dual, and the others follow since
\[
\G(N) = \GL(N)\cap\GR(N) = \GL(M)\cap\GR(M) = \G(M), 
\]
with a similar calculation for $\GLR(N)=\GL(N)\vee\GR(N)$.
\epf

The next statement concerns compositions with $\GLR$, but we note that it says nothing about $\X\circ\GLR$ for $\X=\E,\F,\FL,\FR$.

\begin{lemma}\label{l:G}
For $\X\in\Fu$ we have
\ben
\item \label{G1} $\X\circ\GLR=\begin{cases}
\X &\text{if $\X=\O,\G,\GL,\GR,\GLR$}\\
\GLR &\text{if $\X=\FLR,\I$,}
\end{cases}$
\item \label{G2}
$\GLR\circ\X = \begin{cases}
\O &\text{if $\X=\O,\E$}\\
\G &\text{if $\X=\G,\GL,\GR,\F,\FL,\FR$}\\
\GLR &\text{if $\X=\GLR,\FLR,\I$.}
\end{cases}
$
\een
\end{lemma}

\pf
\firstpfitem{\ref{G1}}  This is clear for $\X=\O$ or $\I$.  For $\X=\G,\GL,\GR,\GLR$ we apply Lemma \ref{l:X} with $N=G_{LR}$:
\[
\X\circ\GLR(M)=\X(G_{LR})=\X(N)=\X(M).
\]
For $\X=\FLR$ we have
\[
G_{LR} \supseteq \FLR(G_{LR}) = \E(G_{LR}) \vee \GLR(G_{LR}) = \E(G_{LR}) \vee G_{LR} = G_{LR},
\]
where we again used Lemma \ref{l:X} in the third step.  Thus, $\FLR(G_{LR})=G_{LR}$: i.e., $\FLR\circ\GLR(M)=\GLR(M)$.

\pfitem{\ref{G2}}  This is again clear for $\X=\O$ or $\I$, and follows from Lemma \ref{l:X} for $\X=\GLR,\FLR$.  For the other choices of $\X$, and writing $X=\X(M)$, the claim follows from previously calculated compositions, in light of $\GLR(X)=\GL(X)\vee\GR(X)$.
\epf

Now we treat compositions with $\FLR$.

\begin{lemma}\label{l:F}
For $\X\in\Fu$ we have
\ben
\item \label{F1} $\X\circ\FLR=\begin{cases}
\FLR &\text{if $\X=\I$}\\
\X &\text{otherwise,}
\end{cases}$
\item \label{F2} $\FLR\circ\X = \begin{cases}
\X &\text{if $\X=\O,\E,\GLR,\FLR$}\\
\G &\text{if $\X=\G,\GL,\GR$}\\
\F &\text{if $\X=\F,\FL,\FR$}\\
\FLR &\text{if $\X=\I$.}
\end{cases}$
\een
\end{lemma}

\pf
\firstpfitem{\ref{F1}}  The $\X=\O,\I$ cases are clear, and Lemma \ref{l:X} (with $N=F_{LR}$) again gives the $\X=\G,\GL,\GR,\GLR$ cases.  The $\X=\E$ case is clear since $E(M)\sub F_{LR}$.  The $\X=\F_\hs$ case follows from the others since
\[
\F_\hs\circ\FLR(M) = \F_\hs(F_{LR}) = \E(F_{LR}) \vee \G_\hs(F_{LR}) = \E(M) \vee \G_\hs(M) = \F_\hs(M).
\]

\pfitem{\ref{F2}}  The $\X=\O,\I$ cases are clear, and the $\X=\GLR$ case is part of Lemma \ref{l:G}\ref{G1}.  The others follow from previously calculated compositions, in light of $\FLR(X)=\E(X)\vee\GLR(X)$, where $X=\X(M)$.
\epf

\subsection{More functors}\label{ss:more}

We have already noted that Table \ref{t:op1} has four missing entries.  At this stage it is conceivable that these missing compositions could be among the functors considered so far, but we will see in Section \ref{ss:size} that they are indeed four new functors.
For now, we simply deal with the missing entries in Table \ref{t:op1} by defining the functors
\begin{equation}\label{e:more}
\Q = \E\circ\GLR \COMMA
\P = \F\circ\GLR \COMMA
\PL = \FL\circ\GLR \COMMA
\PR = \FR\circ\GLR .
\end{equation}
We also define the enlarged set of functors
\[
\Fup = \Fu\cup\{\Q,\P,\PL,\PR\} = \{ \O, \E, \G, \GL, \GR, \GLR, \F, \FL, \FR, \FLR, \Q, \P, \PL, \PR, \I\}.
\]
For a monoid $M$, the functors in \eqref{e:more} yield (at most) four additional submonoids:
\[
\Q(M) = \E(\GLR(M)) \COMMA
\P(M) = \F(\GLR(M)) \COMMA
\PL(M) = \FL(\GLR(M)) \COMMA
\PR(M) = \FR(\GLR(M)) .
\]
Accordingly, we also define
\[
\Lp(M) = \set{\X(M)}{\X\in\Fup}.
\]
We will show in Section \ref{s:enhanced} (see Proposition \ref{p:L2}) that $\Lp(M)$ is a lattice.  Figure \ref{f:lattice2} displays the generic shape of $\Lp(M)$, with the new submonoids shown in red; cf.~Figure \ref{f:lattice}.  The inclusion relations claimed in Figure~\ref{f:lattice2} are all easily verified.  For example, 
\[
\O(M) \sub \Q(M) = \E(\GLR(M)) \sub
\begin{cases}
\E(M) \\
\F(\GLR(M)) = \P(M),
\end{cases}
\]
and for $\hs\not=LR$,
\[
\G_\hs(M) = \G_\hs(\GLR(M)) \sub \F_\hs(\GLR(M)) \sub \F_\hs(M) \implies \G_\hs(M) \sub \P_\hs(M) \sub \F_\hs(M).
\]

\begin{figure}[ht]
\begin{center}
\begin{tikzpicture}[scale=.8]
\node[rounded corners,rectangle,draw,fill=blue!20] (M) at (0,12) {};
\node[rounded corners,rectangle,draw,fill=blue!20] (FLR) at (0,10) {};
\node[rounded corners,rectangle,draw,fill=blue!20] (FL) at (-3,8) {};
\node[rounded corners,rectangle,draw,fill=blue!20] (GLR) at (0,8) {};
\node[rounded corners,rectangle,draw,fill=blue!20] (FR) at (3,8) {};
\node[rounded corners,rectangle,draw,fill=red!20] (PL) at (-3,6) {};
\node[rounded corners,rectangle,draw,fill=blue!20] (F) at (0,6) {};
\node[rounded corners,rectangle,draw,fill=red!20] (PR) at (3,6) {};
\node[rounded corners,rectangle,draw,fill=blue!20] (GL) at (-3,4) {};
\node[rounded corners,rectangle,draw,fill=red!20] (P) at (0,4) {};
\node[rounded corners,rectangle,draw,fill=blue!20] (GR) at (3,4) {};
\node[rounded corners,rectangle,draw,fill=blue!20] (G) at (0,2) {};
\node[rounded corners,rectangle,draw,fill=red!20] (Q) at (6,2) {};
\node[rounded corners,rectangle,draw,fill=blue!20] (E) at (6,4) {};
\node[rounded corners,rectangle,draw,fill=blue!20] (O) at (6,0) {};
\draw (G)--(GR)--(PR)--(GLR)--(PL)--(GL)--(G);
\draw (P)--(PR) (FR)--(FLR)--(FL) (PL)--(P);
\draw[line width=1.5mm,white] (F)--(E) (P)--(Q) (G)--(O);
\draw (O)--(Q) (P)--(G);
\draw (F)--(E) (P)--(Q) (G)--(O);
\draw (FLR)--(GLR) (FL)--(PL) (FR)--(PR) (F)--(P)  (E)--(Q);
\draw[line width=1.5mm,white] (FL)--(F)--(FR);
\draw (FL)--(F)--(FR);
\draw (FLR)--(M);
\node () at (7.5,0) {$=\{1\}$};
\node () at (1.35,12) {$=M$};
\node[rounded corners,rectangle,draw,fill=blue!20] (M) at (0,12) {\small $\I(M)$};
\node[rounded corners,rectangle,draw,fill=blue!20] (FLR) at (0,10) {\small $\F_{LR}(M)$};
\node[rounded corners,rectangle,draw,fill=blue!20] (FL) at (-3,8) {\small $\F_L(M)$};
\node[rounded corners,rectangle,draw,fill=blue!20] (GLR) at (0,8) {\small $\G_{LR}(M)$};
\node[rounded corners,rectangle,draw,fill=blue!20] (FR) at (3,8) {\small $\F_R(M)$};
\node[rounded corners,rectangle,draw,fill=red!20] (PL) at (-3,6) {\small $\P_L(M)$};
\node[rounded corners,rectangle,draw,fill=blue!20] (F) at (0,6) {\small $\F(M)$};
\node[rounded corners,rectangle,draw,fill=red!20] (PR) at (3,6) {\small $\P_R(M)$};
\node[rounded corners,rectangle,draw,fill=blue!20] (GL) at (-3,4) {\small $\G_L(M)$};
\node[rounded corners,rectangle,draw,fill=red!20] (P) at (0,4) {\small $\P(M)$};
\node[rounded corners,rectangle,draw,fill=blue!20] (GR) at (3,4) {\small $\G_R(M)$};
\node[rounded corners,rectangle,draw,fill=blue!20] (G) at (0,2) {\small $\G(M)$};
\node[rounded corners,rectangle,draw,fill=red!20] (Q) at (6,2) {\small $\Q(M)$};
\node[rounded corners,rectangle,draw,fill=blue!20] (E) at (6,4) {\small $\E(M)$};
\node[rounded corners,rectangle,draw,fill=blue!20] (O) at (6,0) {\small $\O(M)$};
\end{tikzpicture}
\caption{The generic shape of the lattice $\Lp(M)$.  In general these submonoids need not be distinct.}
\label{f:lattice2}
\end{center}
\end{figure}

\subsection{More compositions}\label{ss:more2}

Now that we have enlarged our list of functors to $\Fup$, we have a number of further compositions to calculate, namely those of the form $\X\circ\Y$ and $\Y\circ\X$ for $\X\in\Fup$ and $\Y\in\Fup\sm\Fu=\{\Q,\P,\PL,\PR\}$.  These compositions are shown in blue in Table \ref{t:op2}.  All of these entries can be readily verified using Table \ref{t:op1}, associativity of functor composition, and the definition of the new functors.  For example,
\[
\E\circ\Q=\E\circ\E\circ\GLR=\E\circ\GLR=\Q \AND 
\Q\circ\E = \E\circ\GLR\circ\E=\E\circ\O=\O.
\]
As before, some calculations can be performed simultaneously; for example, 
\[
\F_\hs\circ\P_\ds = \F_\hs \circ \F_\ds \circ\GLR = \F\circ\GLR = \P \ANd \P_\hs\circ\P_\ds = \F_\hs\circ\GLR\circ\F_\ds\circ\GLR = \F_\hs\circ\G\circ\GLR = \G\circ\GLR=\G.
\]

Since $\Fup$ is closed under composition (cf.~Table \ref{t:op2}), it is therefore a semigroup, indeed a monoid with identity $\I$.  Note that $\Fup\sm\{\I\}$ is also a semigroup, although it is not a monoid; however, $\FLR$ is a right (but not left) identity element of this subsemigroup.  We will say more about the size and structure of the monoid~$\Fup$ in Sections \ref{ss:size} and \ref{ss:structure}.

\begin{table}[ht]
\begin{center}
\begin{tabular}{l|lllllllllllllll}
\hspace{.5mm}$\circ$ & $\O$ & $\E$ & $\G$  & $\GL$  & $\GR$  & $\GLR$ & $\F$  & $\FL$  & $\FR$  & $\FLR$  & $\Q$ & $\P$ & $\PL$ & $\PR$ & $\I$ \\
\hline
$\O$ & $\O$ & $\O$ & $\O$ & $\O$ & $\O$ & $\O$ & $\O$ & $\O$ & $\O$ & $\O$ & $\O$ & $\O$ & $\O$ & $\O$ & $\O$ \\
$\E$ & $\O$ & $\E$ & $\O$ & $\O$ & $\O$ & {\red $\Q$} & $\E$ & $\E$ & $\E$ & $\E$ & \textcolor{cyan}{ $\Q$} & \textcolor{cyan}{ $\Q$} & \textcolor{cyan}{ $\Q$} & \textcolor{cyan}{ $\Q$} & $\E$ \\
$\G$ & $\O$ & $\O$ & $\G$ & $\G$ & $\G$ & $\G$ & $\G$ & $\G$ & $\G$ & $\G$ & \textcolor{cyan}{ $\O$} & \textcolor{cyan}{ $\G$} & \textcolor{cyan}{ $\G$} & \textcolor{cyan}{ $\G$} & $\G$ \\
$\GL$ & $\O$ & $\O$ & $\G$ & $\G$ & $\G$ & $\GL$ & $\G$ & $\G$ & $\G$ & $\GL$ & \textcolor{cyan}{ $\O$} & \textcolor{cyan}{ $\G$} & \textcolor{cyan}{ $\G$} & \textcolor{cyan}{ $\G$} & $\GL$ \\
$\GR$ & $\O$ & $\O$ & $\G$ & $\G$ & $\G$ & $\GR$ & $\G$ & $\G$ & $\G$ & $\GR$ & \textcolor{cyan}{ $\O$} & \textcolor{cyan}{ $\G$} & \textcolor{cyan}{ $\G$} & \textcolor{cyan}{ $\G$} & $\GR$ \\
$\GLR$ & $\O$ & $\O$ & $\G$ & $\G$ & $\G$ & $\GLR$ & $\G$ & $\G$ & $\G$ & $\GLR$ & \textcolor{cyan}{ $\O$} & \textcolor{cyan}{ $\G$} & \textcolor{cyan}{ $\G$} & \textcolor{cyan}{ $\G$} & $\GLR$ \\
$\F$ & $\O$ & $\E$ & $\G$ & $\G$ & $\G$ & {\red $\P$} & $\F$ & $\F$ & $\F$ & $\F$ & \textcolor{cyan}{ $\Q$} & \textcolor{cyan}{ $\P$} & \textcolor{cyan}{ $\P$} & \textcolor{cyan}{ $\P$} & $\F$ \\
$\FL$ & $\O$ & $\E$ & $\G$ & $\G$ & $\G$ & {\red $\PL$} & $\F$ & $\F$ & $\F$ & $\FL$ & \textcolor{cyan}{ $\Q$} & \textcolor{cyan}{ $\P$} & \textcolor{cyan}{ $\P$} & \textcolor{cyan}{ $\P$} & $\FL$ \\
$\FR$ & $\O$ & $\E$ & $\G$ & $\G$ & $\G$ & {\red $\PR$} & $\F$ & $\F$ & $\F$ & $\FR$ & \textcolor{cyan}{ $\Q$} & \textcolor{cyan}{ $\P$} & \textcolor{cyan}{ $\P$} & \textcolor{cyan}{ $\P$} & $\FR$ \\
$\FLR$ & $\O$ & $\E$ & $\G$ & $\G$ & $\G$ & $\GLR$ & $\F$ & $\F$ & $\F$ & $\FLR$ & \textcolor{cyan}{ $\Q$} & \textcolor{cyan}{ $\P$} & \textcolor{cyan}{ $\P$} & \textcolor{cyan}{ $\P$} & $\FLR$ \\
$\Q$ &$\O$ & \textcolor{cyan}{ $\O$} & \textcolor{cyan}{ $\O$} & \textcolor{cyan}{ $\O$} & \textcolor{cyan}{ $\O$} & \textcolor{cyan}{ $\Q$} & \textcolor{cyan}{ $\O$} & \textcolor{cyan}{ $\O$} & \textcolor{cyan}{ $\O$} & \textcolor{cyan}{ $\Q$} & \textcolor{cyan}{ $\O$} & \textcolor{cyan}{ $\O$} & \textcolor{cyan}{ $\O$} & \textcolor{cyan}{ $\O$} & $\Q$ \\
$\P$ & $\O$ & \textcolor{cyan}{ $\O$} & \textcolor{cyan}{ $\G$} & \textcolor{cyan}{ $\G$} & \textcolor{cyan}{ $\G$} & \textcolor{cyan}{ $\P$} & \textcolor{cyan}{ $\G$} & \textcolor{cyan}{ $\G$} & \textcolor{cyan}{ $\G$} & \textcolor{cyan}{ $\P$} & \textcolor{cyan}{ $\O$} & \textcolor{cyan}{ $\G$} & \textcolor{cyan}{ $\G$} & \textcolor{cyan}{ $\G$} & $\P$ \\
$\PL$ & $\O$ & \textcolor{cyan}{ $\O$} & \textcolor{cyan}{ $\G$} & \textcolor{cyan}{ $\G$} & \textcolor{cyan}{ $\G$} & \textcolor{cyan}{ $\PL$} & \textcolor{cyan}{ $\G$} & \textcolor{cyan}{ $\G$} & \textcolor{cyan}{ $\G$} & \textcolor{cyan}{ $\PL$} & \textcolor{cyan}{ $\O$} & \textcolor{cyan}{ $\G$} & \textcolor{cyan}{ $\G$} & \textcolor{cyan}{ $\G$} & $\PL$ \\
$\PR$ & $\O$ & \textcolor{cyan}{ $\O$} & \textcolor{cyan}{ $\G$} & \textcolor{cyan}{ $\G$} & \textcolor{cyan}{ $\G$} & \textcolor{cyan}{ $\PR$} & \textcolor{cyan}{ $\G$} & \textcolor{cyan}{ $\G$} & \textcolor{cyan}{ $\G$} & \textcolor{cyan}{ $\PR$} & \textcolor{cyan}{ $\O$} & \textcolor{cyan}{ $\G$} & \textcolor{cyan}{ $\G$} & \textcolor{cyan}{ $\G$} & $\PR$ \\
$\I$ &  $\O$ & $\E$ & $\G$  & $\GL$  & $\GR$  & $\GLR$ & $\F$  & $\FL$  & $\FR$  & $\FLR$  & $\Q$ & $\P$ & $\PL$ & $\PR$ & $\I$ \\
\end{tabular}
\caption{Composition of the functors from $\Fup$.}
\label{t:op2}
\end{center}
\end{table}

\subsection{Size}\label{ss:size}

We now know that the set $\Fup$ is a monoid under composition, and that its size is at most $15$.  We also know that $|\Fup|\geq11$, since $|\Lp(M)|\geq|\L(M)|=11$ for  $M$ of type $(0,0,0,0)$; cf.~Figure \ref{f:L2}.
To show that the size of $\Fup$ is in fact $15$, as we will in Proposition \ref{p:15} below, we will construct a monoid~$M$ such that~$\Lp(M)$ has size $15$.  We begin by showing that the functors from $\Fup$ respect the direct product operation:

\begin{lemma}\label{l:DP2}
For any $\X\in\Fup$, and for any two monoids $M$ and $N$, we have
\[
\X(M\times N) = \X(M) \times \X(N).
\]
\end{lemma}

\pf
In light of Lemma \ref{l:DP}, it suffices to demonstrate this for any $\X\in\Fup\sm\Fu$.  For any such $\X$, we have $\X=\Y\circ\Z$ for some $\Y,\Z\in\Fu$.  Two applications of Lemma \ref{l:DP} then give
\[
\X(M\times N) = \Y(\Z(M\times N)) = \Y(\Z(M)\times\Z(N)) = \Y(\Z(M)) \times \Y(\Z(N)) = \X(M) \times \X(N). \qedhere
\]
\epf

Table \ref{t:GNB_small} listed the submonoids $\X(M)$, $\X\in\Fu$, for various monoids $M$ defined in Section \ref{ss:egs}.   Table~\ref{t:GNB_smaller} gives the submonoids $\X(M)$ for the additional functors $\X\in\Fup\sm\Fu$.
The entries for $M=G$, $E$ and $\Pos$ are clear, while those for $M=B$ and $B^0$ follow quickly from Table \ref{t:GNB_small} and the fact that $\GLR(B^0)=\GLR(B)=B$.

\begin{table}[ht]
\begin{center}
\begin{tabular}{l|ccccc}
$\X$ & $\X(G)$ & $\X(E)$ & $\X(\Pos)$ & $\X(B)$  & $\X(B^0)$ \\
\hline
$\Q$ & $\{1\}$ & $\{1\}$& $\{1\}$ & $\set{a^mb^m}{m\geq0}$  & $\set{a^mb^m}{m\geq0}$ \\
$\P$ & $G$ & $\{1\}$& $\{1\}$ & $\set{a^mb^m}{m\geq0}$  & $\set{a^mb^m}{m\geq0}$ \\
$\PL$ & $G$ & $\{1\}$& $\{1\}$ & $\set{a^mb^n}{m\geq n}$  & $\set{a^mb^n}{m\geq n}$ \\
$\PR$ & $G$ & $\{1\}$& $\{1\}$ & $\set{a^mb^n}{m\leq n}$  & $\set{a^mb^n}{m\leq n}$ \\
\end{tabular}
\caption{The submonoids $\X(M)$, $\X\in\Fup\sm\Fu$, for $M=G$ (a group), $M=E$ (an idempotent-generated monoid), $M=\Pos$ (the positive integers under multiplication), $M=B$ (the bicyclic monoid) and $M=B^0$ (the bicyclic monoid with a zero adjoined); cf.~Table \ref{t:GNB_small}.}
\label{t:GNB_smaller}
\end{center}
\end{table}

\begin{prop}\label{p:15}
The monoid $\Fup$ has size $15$.
\end{prop}

\pf
Consider the monoid $M = G\times E\times\Pos\times B$, where $G$ is a nontrivial group, $E$ a nontrivial idempotent-generated monoid, $\Pos$ the positive integers under multiplication, and $B$ the bicyclic monoid.  By consulting Tables \ref{t:GNB_small} and \ref{t:GNB_smaller}, and keeping in mind that $\X(M) = \X(G)\times\X(E)\times\X(\Pos)\times\X(B)$ for all $\X\in\Fup$ (cf.~Lemma~\ref{l:DP2}), one may easily check that $\Lp(M)$ has size $15$.
\epf

\subsection{Structure}\label{ss:structure}

Now that we know the size of the monoid $\Fup$, it is natural to seek more information about its algebraic structure.  The most common way to structurally decompose a semigroup is by using Green's relations.  These were defined in Section \ref{ss:Green}, but we recall a number of additional definitions here.  

For a monoid $M$, Green's $\gJ$-preorder is the relation $\leqJ$ defined by $x\leqJ y \iff MxM\sub MyM$.  Again, this may be reformulated in terms of divisibility: $x\leqJ y\iff x=ayb$ for some $a,b\in M$.  Green's $\gJ$ relation (as defined in Section \ref{ss:Green}) is then given by ${\gJ}={\leqJ}\cap{\geqJ}$.  Note that $\leqJ$ is a partial order if and only if $M$ is $\gJ$-trivial.

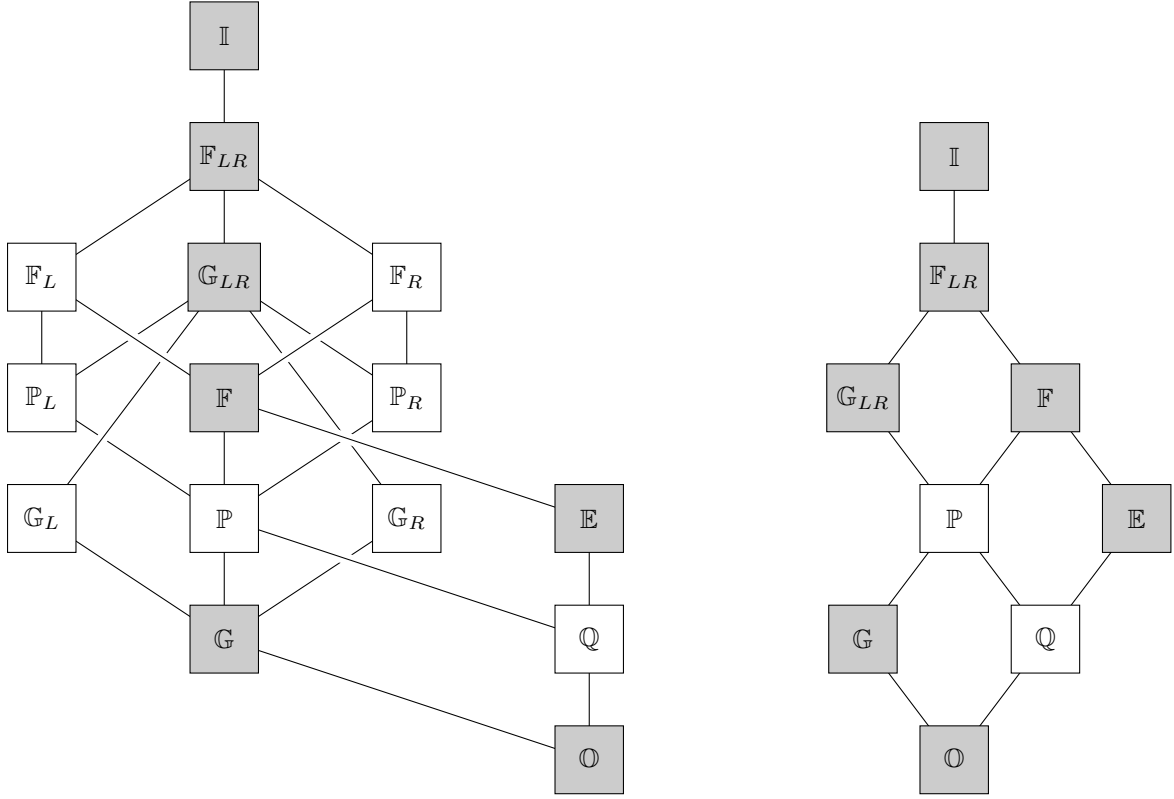
\begin{figure}[ht]
\begin{center}
\begin{tikzpicture}[scale=0.8,
block/.style={
draw,
fill=white,
rectangle, 
minimum width=.9cm,
minimum height=.9cm,
font=\small}]
\nc\hh{1}
\begin{scope}
\node (I) at (0,12*\hh) {};
\node (FLR) at (0,10*\hh) {};
\node (FL) at (-3,8*\hh) {};
\node (GLR) at (0,8*\hh) {};
\node (FR) at (3,8*\hh) {};
\node (PL) at (-3,6*\hh) {};
\node (F) at (0,6*\hh) {};
\node (PR) at (3,6*\hh) {};
\node (GL) at (-3,4*\hh) {};
\node (P) at (0,4*\hh) {};
\node (GR) at (3,4*\hh) {};
\node (E) at (6,4*\hh) {};
\node (G) at (0,2*\hh) {};
\node (Q) at (6,2*\hh) {};
\node (O) at (6,0) {};
\draw (O)--(G)--(P)--(PL)--(GLR)--(PR)--(P);
\draw[line width=1.5mm,white] (G)--(GR)--(GLR)--(GL)--(G);
\draw (G)--(GR)--(GLR)--(GL)--(G);
\draw[line width=1.5mm,white] (Q)--(P)--(F) ;
\draw (Q)--(P)--(F);
\draw (O)--(Q)--(E)--(F);
\draw (PL)--(FL)--(FLR)--(FR)--(PR);
\draw (GLR)--(FLR)--(I);
\draw[line width=1.5mm,white] (FL)--(F)--(FR) (E)--(F);
\draw (FL)--(F)--(FR) (E)--(F);
\node[block,fill=black!20] (I) at (0,12*\hh) {$\I$};
\node[block,fill=black!20] (FLR) at (0,10*\hh) {$\FLR$};
\node[block] (FL) at (-3,8*\hh) {$\FL$};
\node[block,fill=black!20] (GLR) at (0,8*\hh) {$\GLR$};
\node[block] (FR) at (3,8*\hh) {$\FR$};
\node[block] (PL) at (-3,6*\hh) {$\PL$};
\node[block,fill=black!20] (F) at (0,6*\hh) {$\F$};
\node[block] (PR) at (3,6*\hh) {$\PR$};
\node[block] (GL) at (-3,4*\hh) {$\GL$};
\node[block] (P) at (0,4*\hh) {$\P$};
\node[block] (GR) at (3,4*\hh) {$\GR$};
\node[block,fill=black!20] (E) at (6,4*\hh) {$\E$};
\node[block,fill=black!20] (G) at (0,2*\hh) {$\G$};
\node[block] (Q) at (6,2*\hh) {$\Q$};
\node[block,fill=black!20] (O) at (6,0) {$\O$};
\end{scope}
\begin{scope}[shift={(12,0)}]
\node (I) at (0,10*\hh) {};
\node (FLR) at (0,8*\hh) {};
\node (GLR) at (-1.5,6*\hh) {};
\node (F) at (1.5,6*\hh) {};
\node (P) at (0,4*\hh) {};
\node (E) at (3,4*\hh) {};
\node (G) at (-1.5,2*\hh) {};
\node (Q) at (1.5,2*\hh) {};
\node (O) at (0,0) {};
\draw(O)--(G)--(P)--(GLR)--(FLR)--(I) (O)--(Q)--(E)--(F) (Q)--(P)--(F)--(FLR);
\node[block,fill=black!20] (I) at (0,10*\hh) {$\I$};
\node[block,fill=black!20] (FLR) at (0,8*\hh) {$\FLR$};
\node[block,fill=black!20] (GLR) at (-1.5,6*\hh) {$\GLR$};
\node[block,fill=black!20] (F) at (1.5,6*\hh) {$\F$};
\node[block] (P) at (0,4*\hh) {$\P$};
\node[block,fill=black!20] (E) at (3,4*\hh) {$\E$};
\node[block,fill=black!20] (G) at (-1.5,2*\hh) {$\G$};
\node[block] (Q) at (1.5,2*\hh) {$\Q$};
\node[block,fill=black!20] (O) at (0,0) {$\O$};
\end{scope}
\end{tikzpicture}
\caption{The divisibility order in the monoids $\Fup$ and $\E(\Fup)$: left and right, respectively.}
\label{f:div}
\end{center}
\end{figure}

Using the composition table (cf.~Table \ref{t:op2}), the computational algebra system GAP \cite{GAP} can perform many calculations in the monoid $\Fup$.  For example, GAP verifies that the monoid $\Fup$ is in fact $\gJ$-trivial, and hence $\gL$-, $\gR$-, $\gH$- and $\gD$-trivial as well.  GAP was also used to produce Figure \ref{f:div}, which displays the $\leqJ$ order in $\Fup$.  In fact, since $\Fup$ is $\gJ$-trivial, Figure~\ref{f:div} is also the so-called \emph{eggbox diagram} of $\Fup$, as defined for example in \cite[Section 1.2]{Hig}.  

As is customary, the idempotents of $\Fup$ (i.e., the functors $\X=\X\circ\X$) are coloured grey in Figure \ref{f:div}.  The idempotent-generated submonoid 
\begin{equation}\label{e:EF}
\E(\Fup) = \{ \O, \E, \G, \GLR, \F, \FLR, \Q, \P, \I\}
\end{equation}
is also pictured in Figure \ref{f:div}, along with its $\leqJ$ ordering, again with assistance from GAP.

GAP also shows that $\Fup$ has 2904 subsemigroups (exactly half of which are submonoids), and 1613 congruences, of which 76 are principal.  (A congruence on a semigroup is an equivalence relation compatible with the product; these are used to form quotient semigroups; see \cite[Section 1.5]{Howie} for more details.)

It is interesting to compare the two posets $(\Lp(M),{\sub})$ and $(\Fup,{\leqJ})$, which are pictured in Figures \ref{f:lattice2} and \ref{f:div}, respectively (the former in the generic case).  Although there are certainly some superficial similarities between them, the two posets are not (quite) isomorphic.  For example, we have $\GL(M)\sub\PL(M)$ for any monoid~$M$ (cf.~Figure \ref{f:lattice2}), while $\GL\not\leqJ\PL$ in~$\Fup$ (cf.~Figure~\ref{f:div}).  We can also see that $\GL\not\leqJ\PL$ directly; using Table~\ref{t:op2}, it is easy to verify that for any $\X,\Y\in\Fup$, we have $\X\circ\PL\circ\Y\in\{\O,\G,\Q,\P,\PL\}$, which means that $\GL\not\in\Fup\circ\PL\circ\Fup$.

\subsection{The lattice of the monoid of functors}\label{ss:LF}

Since $\Fup$ is a monoid, it is natural to calculate its associated lattice $\L(\Fup)$.  Consulting Table \ref{t:op2}, we see that the only solution in $\Fup$ to $\X\circ\Y=\I$ is $\X=\Y=\I$, so it follows that $\G(\Fup)=\GL(\Fup)=\{\I\}$: i.e., that $\Tone(\Fup)=\Tfour(\Fup)=1$.  From \eqref{e:EF} we have $\{\I\}\subsetneq\E(\Fup)=\F(\Fup)\subsetneq\Fup$; note that $\E(\Fup)=\F(\Fup)$ because $\G(\Fup)=\{\I\}$; cf.~Lemma \ref{l:EG}\ref{EG2}.  It follows that $\Ttwod(\Fup)=\Tthreed(\Fup)=0$.  All of the above shows that the monoid $\Fup$ has type $\T(\Fup) = (1,0,0,1)$, and so
\[
\L(\Fup) = \big\{ \{\I\}, \E(\Fup), \Fup \big\}
\]
is the three-element chain displayed in the second diagram on the top row of Figure \ref{f:L1}.

\section{An enhanced lattice invariant?}\label{s:enhanced}

Section \ref{s:class} concerned the lattice $\L(M)=\set{\X(M)}{\X\in\Fu}$ consisting of the submonoids of a monoid $M$ arising from the functors from $\Fu$.  Section \ref{s:F} concerned the monoid $\Fup$ of functors generated by $\Fu$, and we defined $\Lp(M)=\set{\X(M)}{\X\in\Fup}$.  As promised earlier, we now show that $\Lp(M)$ is a lattice.

\begin{prop}\label{p:L2}
For any monoid $M$, the set $\Lp(M)$ is a finite $\vee$-subsemilattice of $\Sub(M)$, with top element~$\I(M)=M$ and bottom element~$\O(M)=\{1\}$.  Consequently, $\Lp(M)$ is a lattice.
\end{prop}

\pf
As in the proof of Proposition \ref{p:L}, it suffices to show that $\Lp(M)$ is closed under $\vee$.  To do so, let~$\X,\Y\in\Fup$; we must show that 
\begin{equation}\label{e:XY}
\X(M)\vee\Y(M) = \Z(M) \qquad\text{for some $\Z\in\Fup$.}
\end{equation}
In light of Proposition \ref{p:L}, and by commutativity of $\vee$, we may assume that $\X\in\Fup\sm\Fu=\{\Q,\P,\PL,\PR\}$.

If $\Y=\I$ or $\O$, then \eqref{e:XY} is clear; we take $\Z=\I$ or $\X$, respectively.

Next suppose $\Y$ is one of $\G_\hs$, $\Q$ or $\P_\ds$.  Then, consulting Table \ref{t:op2}, we see that $\X=\U\circ\GLR$ and $\Y=\V\circ\GLR$ for some $\U,\V\in\Fu$.  But then
\[
\X(M) \vee \Y(M) = \U(G_{LR}) \vee \V(G_{LR}) = \W(G_{LR}) = \W\circ\GLR(M) \qquad\text{for some $\W\in\Fu$,}
\]
using Proposition \ref{p:L} (applied in the monoid $G_{LR}$) in the second step.  We then take $\Z=\W\circ\GLR\in\Fup$.

Next suppose $\Y=\E$.  If $\X=\Q$, then \eqref{e:XY} is clear since $Q\sub E$.  Now suppose $\X=\P_\hs$.  Then
\[
P_\hs \vee E = \F_\hs(G_{LR}) \vee E = \E(G_{LR}) \vee \G_\hs(G_{LR}) \vee E = G_\hs\vee E = F_\hs,
\]
where we used $\E(G_{LR})\sub E$ and $\G_\hs(G_{LR})=G_\hs$ in the third step.  Thus, we may take $\Z=\F_\hs$ in this case.

Finally, suppose $\Y=\F_\hs$.  Again \eqref{e:XY} is clear for $\X=\Q$, as $Q\sub F_\hs$, so suppose $\X=\P_\ds$.  Then writing ${G_\ds\vee G_\hs=G_\ss}$, we have
\[
P_\ds\vee F_\hs = \F_\ds(G_{LR}) \vee F_\hs = \E(G_{LR}) \vee \G_\ds(G_{LR}) \vee E \vee G_\hs = G_\ds\vee E \vee G_\hs = E\vee G_\ss = F_\ss,
\]
so we may take $\Z=\F_\ss$ in this case.
\epf

We now have two lattice invariants $\L(M)$ and $\Lp(M)$, associated to a monoid $M$.  Given that $\Lp(M)$ is defined in terms of a larger set of functors, one might hope that it allows us to distinguish monoids not distinguished by $\L(M)$.  However, it follows from the results of this section that this is not the case.  In Section~\ref{ss:more3} we prove some preliminary results about collapse in the enhanced lattice $\Lp(M)$, and then we classify the lattices $\Lp(M)$ in Section \ref{ss:class2}.

\subsection{More collapse}\label{ss:more3}

We begin with some results analogous to Lemmas \ref{l:GLR} and \ref{l:EG}, but involving the functors from $\Fup\sm\Fu$.  For a monoid $M$, we continue to use the abbreviations \eqref{e:G} and \eqref{e:F}, as well as
\[
Q=\Q(M) \COMMA P=\P(M) \COMMA P_L=\PL(M) \COMMA P_R=\PR(M).
\]

\begin{lemma}\label{l:PGLR}
For a monoid $M$, conditions \ref{GLR1}--\ref{GLR5} of Lemma \ref{l:GLR} are also equivalent to each of the following:
\ben\addtocounter{enumi}{7}
\item \label{PGLR1} $P$, $P_L$, $P_R$ and $G_{LR}$ are not all equal,
\item \label{PGLR2} $P$, $P_L$, $P_R$ and $G_{LR}$ are pairwise distinct.
\een
\end{lemma}

\pf
Writing $N=\GLR(M)$, note that
\[
P = \F(N) \COMMA P_L=\FL(N) \COMMA P_R=\FR(N) \COMMA G_{LR}=\FLR(N).
\]
Thus, the equivalence of \ref{PGLR1} and \ref{PGLR2} follows from the equivalence of \ref{FLR1} and \ref{FLR2} in the monoid $N$.

\pfitem{\ref{GLR2}$\implies$\ref{PGLR1}}  Aiming to prove the contrapositive, suppose \ref{PGLR1} does not hold.  In particular, we have $P=G_{LR}$: i.e., $\P(M)=\GLR(M)$.  But then
\[
G = \G(M) = \GLR(\P(M)) = \GLR(\GLR(M)) = \GLR(M) = G_{LR},
\]
so that \ref{GLR2} does not hold.

\pfitem{\ref{PGLR2}$\implies$\ref{GLR1}}  Again we prove the contrapositive.  If \ref{GLR1} does not hold, then $G=G_{LR}$, from which it follows that
\[
P = \P(M) = \F(\GLR(M)) = \F(\G(M)) = \G(M) = G = G_{LR},
\]
so that \ref{PGLR2} does not hold.
\epf

\begin{lemma}\label{l:PQ}
For any monoid $M$ we have
\ben
\item \label{PQ1} $G=\{1\} \iff P=Q \iff F=E$,
\item \label{PQ2} $E=Q \iff F=P \iff F_L=P_L \iff F_R=P_R \iff F_{LR}=G_{LR}$.
\een
\end{lemma}

\pf
\firstpfitem{\ref{PQ1}}
In light of Lemma \ref{l:EG}\ref{EG2}, it is enough to show that $G=\{1\} \iff P=Q$.
If $G=\{1\}$, then 
\[
P = \F(G_{LR}) = \E(G_{LR}) \vee \G(G_{LR}) = Q\vee G = Q\vee\{1\} = Q.
\]
Conversely, if $P=Q$, then
\[
G = \G(M) = \G\circ\P(M) = \G\circ\Q(M) = \O(M) = \{1\}.  
\]

\pfitem{\ref{PQ2}}
For convenience during this part of the proof, we will write $\PLR=\GLR$ and $P_{LR}=G_{LR}$.  So we wish to show that $E=Q \iff F_\hs=P_\hs$ for any subscript $\hs$.  First, if $E=Q$, then 
\[
F_\hs = E\vee G_\hs = Q\vee G_\hs = \E(G_{LR})\vee\G_\hs(G_{LR}) = \F_\hs(G_{LR}) = P_\hs.
\]
(Note that the last step holds by definition apart from the $\hs=LR$ case, when it follows instead from Lemma~\ref{l:G}\ref{G1} and the $P_{LR}=G_{LR}$ convention.)  Conversely, if $F_\hs=P_\hs$ for some choice of $\hs$, then
\[
E = \E(M) = \E\circ\F_\hs(M) = \E\circ\P_\hs(M) = \Q(M) = Q.  \qedhere
\]
\epf

\subsection{Classification of enhanced lattice invariants}\label{ss:class2}

We now wish to classify the enhanced lattice invariants $\Lp(M)$, for monoids $M$.  To do so, we will again use the type of $M$, as defined in Section \ref{s:class}.  

First note that if $\Tone(M)=1$, then $G=G_L$ and so $G_{LR}=G$ (cf.~Lemma \ref{l:GLR}), so it follows that
\[
Q = \E(G_{LR}) = \E(G) = \{1\} \AND P_\hs=\F_\hs(G_{LR}) = \F_\hs(G) = G,
\]
which means that $\Lp(M)=\L(M)$ in this case.

If $\Tthree(M)=1$, then $G_{LR}=F_{LR}$, and this time
\[
Q = 
\E(G_{LR}) = \E(F_{LR}) = E 
\AND 
P_\hs = 
\F_\hs(G_{LR}) = \F_\hs(F_{LR}) = F_\hs,
\]
so that $\Lp(M)=\L(M)$ in this case as well.

This leaves us to consider monoids $M$ of type $\T(M)=(0,i,0,j)$.  As explained in Section~\ref{s:class}, the seven submonoids of $M$ listed in \eqref{e:7} are distinct (as $\Tone(M)=0$).  By Lemma \ref{l:PGLR}, the four submonoids $P,P_L,P_R,G_{LR}$ are distinct as well.  Because also $F_{LR}\not= G_{LR}$ (as $\Tthree(M)=0$), it follows from Lemma \ref{l:PQ}\ref{PQ2} that the following containments are strict:
\[
Q\subsetneq E \COMMA
P\subsetneq F \COMMA
P_L\subsetneq F_L \COMMA
P_R\subsetneq F_R \COMMA
G_{LR}\subsetneq F_{LR} .
\]
We claim that the following containments are also strict:
\[
\{1\} \subsetneq Q \COMMA
G \subsetneq P \COMMA
G_L \subsetneq P_L \COMMA
G_R \subsetneq P_R .
\]
Indeed, to see this, note first that $\{1\}$, $G$, $G_L$ and $G_R$ have only one idempotent (cf.~Lemma \ref{l:GE}\ref{GE1}).  On the other hand, $Q=\E(\GLR(M))$ contains infinitely many idempotents (cf.~Lemma \ref{l:GLR}, and note that $G\not=G_L$ since $\Tone(M)=0$), and so too do each of $P$, $P_L$ and $P_R$, since all three of these contain $Q$.  This completes the proof of the claim.  All of the above shows that the following eleven submonoids of $M$ are distinct:
\[
\text{$G$, $G_L$, $G_R$, $G_{LR}$, $P$, $P_L$, $P_R$, $F$, $F_L$, $F_R$ and $F_{LR}$.}
\]
These submonoids are shaded red in Figure \ref{f:lattice_red2}, which gives the generic shape of $\Lp(M)$ in the case that~$\T(M)=(0,i,0,j)$.

\begin{figure}[ht]
\begin{center}
\begin{tikzpicture}[scale=.6]
\node[rounded corners,rectangle,draw,fill=blue!20] (M) at (0,12) {};
\node[rounded corners,rectangle,draw,fill=red!20] (FLR) at (0,10) {};
\node[rounded corners,rectangle,draw,fill=red!20] (FL) at (-3,8) {};
\node[rounded corners,rectangle,draw,fill=red!20] (GLR) at (0,8) {};
\node[rounded corners,rectangle,draw,fill=red!20] (FR) at (3,8) {};
\node[rounded corners,rectangle,draw,fill=red!20] (PL) at (-3,6) {};
\node[rounded corners,rectangle,draw,fill=red!20] (F) at (0,6) {};
\node[rounded corners,rectangle,draw,fill=red!20] (PR) at (3,6) {};
\node[rounded corners,rectangle,draw,fill=red!20] (GL) at (-3,4) {};
\node[rounded corners,rectangle,draw,fill=red!20] (P) at (0,4) {};
\node[rounded corners,rectangle,draw,fill=red!20] (GR) at (3,4) {};
\node[rounded corners,rectangle,draw,fill=red!20] (G) at (0,2) {};
\node[rounded corners,rectangle,draw,fill=blue!20] (E) at (6,4) {};
\node[rounded corners,rectangle,draw,fill=blue!20] (Q) at (6,2) {};
\node[rounded corners,rectangle,draw,fill=blue!20] (1) at (6,0) {};
\draw[ultra thick] (FLR)--(FL)--(PL)--(GL)--(G)--(GR)--(PR)--(FR)--(FLR)--(GLR)--(PL)--(P)--(PR)--(GLR) (F)--(P)--(G) (E)--(Q)--(1);
\draw[line width=2mm,white] (FL)--(F)--(FR);
\draw[ultra thick] (FL)--(F)--(FR);
\draw (M)--(FLR)--(GLR);
\draw[line width=1.5mm,white] (P)--(Q) (F)--(E);
\draw (F)--(E) (P)--(Q) (G)--(1);
\node[rounded corners,rectangle,draw,fill=blue!20] (M) at (0,12) {\small $M$};
\node[rounded corners,rectangle,draw,fill=red!20] (FLR) at (0,10) {\small $F_{LR}$};
\node[rounded corners,rectangle,draw,fill=red!20] (FL) at (-3,8) {\small $F_L$};
\node[rounded corners,rectangle,draw,fill=red!20] (GLR) at (0,8) {\small $G_{LR}$};
\node[rounded corners,rectangle,draw,fill=red!20] (FR) at (3,8) {\small $F_R$};
\node[rounded corners,rectangle,draw,fill=red!20] (PL) at (-3,6) {\small $P_L$};
\node[rounded corners,rectangle,draw,fill=red!20] (F) at (0,6) {\small $F$};
\node[rounded corners,rectangle,draw,fill=red!20] (PR) at (3,6) {\small $P_R$};
\node[rounded corners,rectangle,draw,fill=red!20] (GL) at (-3,4) {\small $G_L$};
\node[rounded corners,rectangle,draw,fill=red!20] (P) at (0,4) {\small $P$};
\node[rounded corners,rectangle,draw,fill=red!20] (GR) at (3,4) {\small $G_R$};
\node[rounded corners,rectangle,draw,fill=red!20] (G) at (0,2) {\small $G$};
\node[rounded corners,rectangle,draw,fill=blue!20] (E) at (6,4) {\small $E$};
\node[rounded corners,rectangle,draw,fill=blue!20] (Q) at (6,2) {\small $Q$};
\node[rounded corners,rectangle,draw,fill=blue!20] (1) at (6,0) {\small $1$};
\end{tikzpicture}
\caption{The generic shape of the lattice $\Lp(M)$ when $M$ has type $\T(M)=(0,i,0,j)$.  The submonoids shaded red are distinct, and thick lines indicate proper containment.}
\label{f:lattice_red2}
\end{center}
\end{figure}

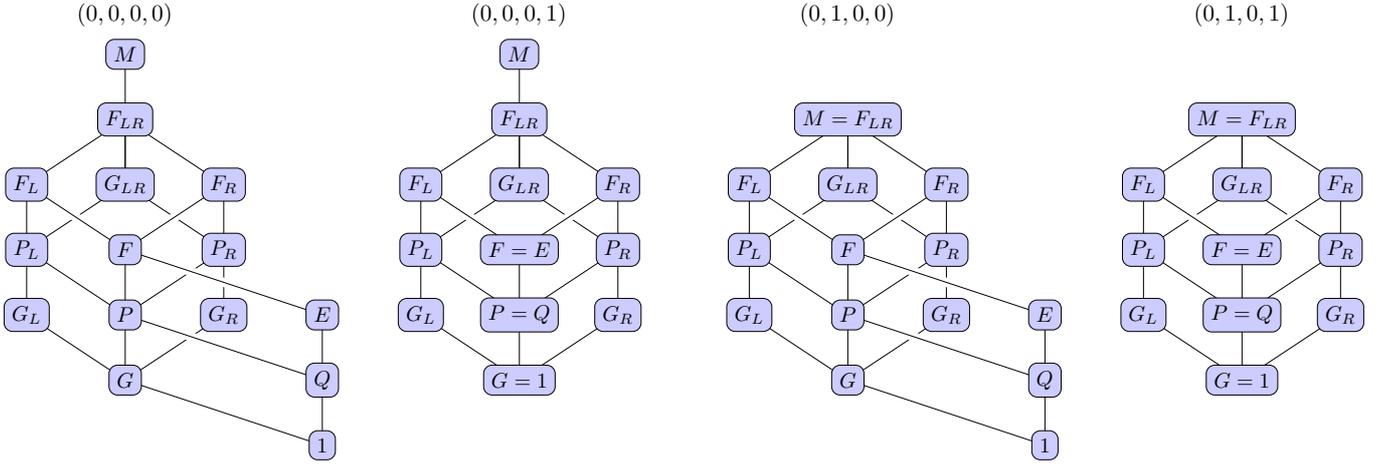
\begin{figure}[ht]
\begin{center}
\scalebox{0.8}{
\begin{tikzpicture}[scale=.54]
\begin{scope}[shift={(0,0)}]
\node[rounded corners,rectangle,draw,fill=blue!20] (M) at (0,12) {};
\node[rounded corners,rectangle,draw,fill=blue!20] (FLR) at (0,10) {};
\node[rounded corners,rectangle,draw,fill=blue!20] (FL) at (-3,8) {};
\node[rounded corners,rectangle,draw,fill=blue!20] (GLR) at (0,8) {};
\node[rounded corners,rectangle,draw,fill=blue!20] (FR) at (3,8) {};
\node[rounded corners,rectangle,draw,fill=blue!20] (PL) at (-3,6) {};
\node[rounded corners,rectangle,draw,fill=blue!20] (F) at (0,6) {};
\node[rounded corners,rectangle,draw,fill=blue!20] (PR) at (3,6) {};
\node[rounded corners,rectangle,draw,fill=blue!20] (GL) at (-3,4) {};
\node[rounded corners,rectangle,draw,fill=blue!20] (P) at (0,4) {};
\node[rounded corners,rectangle,draw,fill=blue!20] (GR) at (3,4) {};
\node[rounded corners,rectangle,draw,fill=blue!20] (G) at (0,2) {};
\node[rounded corners,rectangle,draw,fill=blue!20] (E) at (6,4) {};
\node[rounded corners,rectangle,draw,fill=blue!20] (Q) at (6,2) {};
\node[rounded corners,rectangle,draw,fill=blue!20] (1) at (6,0) {};
\draw (FLR)--(FL)--(PL)--(GL)--(G)--(GR)--(PR)--(FR)--(FLR)--(GLR)--(PL)--(P)--(PR)--(GLR) (F)--(P)--(G) (E)--(Q)--(1);
\draw[line width=1.5mm,white] (FL)--(F)--(FR);
\draw (FL)--(F)--(FR);
\draw (M)--(FLR)--(GLR);
\draw[line width=1.5mm,white] (P)--(Q) (F)--(E);
\draw (F)--(E) (P)--(Q) (G)--(1);
\node[rounded corners,rectangle,draw,fill=blue!20] (M) at (0,12) {\small $M$};
\node[rounded corners,rectangle,draw,fill=blue!20] (FLR) at (0,10) {\small $F_{LR}$};
\node[rounded corners,rectangle,draw,fill=blue!20] (FL) at (-3,8) {\small $F_L$};
\node[rounded corners,rectangle,draw,fill=blue!20] (GLR) at (0,8) {\small $G_{LR}$};
\node[rounded corners,rectangle,draw,fill=blue!20] (FR) at (3,8) {\small $F_R$};
\node[rounded corners,rectangle,draw,fill=blue!20] (PL) at (-3,6) {\small $P_L$};
\node[rounded corners,rectangle,draw,fill=blue!20] (F) at (0,6) {\small $F$};
\node[rounded corners,rectangle,draw,fill=blue!20] (PR) at (3,6) {\small $P_R$};
\node[rounded corners,rectangle,draw,fill=blue!20] (GL) at (-3,4) {\small $G_L$};
\node[rounded corners,rectangle,draw,fill=blue!20] (P) at (0,4) {\small $P$};
\node[rounded corners,rectangle,draw,fill=blue!20] (GR) at (3,4) {\small $G_R$};
\node[rounded corners,rectangle,draw,fill=blue!20] (G) at (0,2) {\small $G$};
\node[rounded corners,rectangle,draw,fill=blue!20] (E) at (6,4) {\small $E$};
\node[rounded corners,rectangle,draw,fill=blue!20] (Q) at (6,2) {\small $Q$};
\node[rounded corners,rectangle,draw,fill=blue!20] (1) at (6,0) {\small $1$};
\node () at (0,13.2) {$(0,0,0,0)$};
\end{scope}
\begin{scope}[shift={(12,0)}]
\node[rounded corners,rectangle,draw,fill=blue!20] (M) at (0,12) {};
\node[rounded corners,rectangle,draw,fill=blue!20] (FLR) at (0,10) {};
\node[rounded corners,rectangle,draw,fill=blue!20] (FL) at (-3,8) {};
\node[rounded corners,rectangle,draw,fill=blue!20] (GLR) at (0,8) {};
\node[rounded corners,rectangle,draw,fill=blue!20] (FR) at (3,8) {};
\node[rounded corners,rectangle,draw,fill=blue!20] (PL) at (-3,6) {};
\node[rounded corners,rectangle,draw,fill=blue!20] (F) at (0,6) {};
\node[rounded corners,rectangle,draw,fill=blue!20] (PR) at (3,6) {};
\node[rounded corners,rectangle,draw,fill=blue!20] (GL) at (-3,4) {};
\node[rounded corners,rectangle,draw,fill=blue!20] (P) at (0,4) {};
\node[rounded corners,rectangle,draw,fill=blue!20] (GR) at (3,4) {};
\node[rounded corners,rectangle,draw,fill=blue!20] (G) at (0,2) {};
\draw (FLR)--(FL)--(PL)--(GL)--(G)--(GR)--(PR)--(FR)--(FLR)--(GLR)--(PL)--(P)--(PR)--(GLR) (F)--(P)--(G);
\draw[line width=1.5mm,white] (FL)--(F)--(FR);
\draw (FL)--(F)--(FR);
\draw (M)--(FLR)--(GLR);
\node[rounded corners,rectangle,draw,fill=blue!20] (M) at (0,12) {\small $M$};
\node[rounded corners,rectangle,draw,fill=blue!20] (FLR) at (0,10) {\small $F_{LR}$};
\node[rounded corners,rectangle,draw,fill=blue!20] (FL) at (-3,8) {\small $F_L$};
\node[rounded corners,rectangle,draw,fill=blue!20] (GLR) at (0,8) {\small $G_{LR}$};
\node[rounded corners,rectangle,draw,fill=blue!20] (FR) at (3,8) {\small $F_R$};
\node[rounded corners,rectangle,draw,fill=blue!20] (PL) at (-3,6) {\small $P_L$};
\node[rounded corners,rectangle,draw,fill=blue!20] (F) at (0,6) {\small $F=E$};
\node[rounded corners,rectangle,draw,fill=blue!20] (PR) at (3,6) {\small $P_R$};
\node[rounded corners,rectangle,draw,fill=blue!20] (GL) at (-3,4) {\small $G_L$};
\node[rounded corners,rectangle,draw,fill=blue!20] (P) at (0,4) {\small $P=Q$};
\node[rounded corners,rectangle,draw,fill=blue!20] (GR) at (3,4) {\small $G_R$};
\node[rounded corners,rectangle,draw,fill=blue!20] (G) at (0,2) {\small $G=1$};
\node () at (0,13.2) {$(0,0,0,1)$};
\end{scope}
\begin{scope}[shift={(22,0)}]
\node[rounded corners,rectangle,draw,fill=blue!20] (FLR) at (0,10) {};
\node[rounded corners,rectangle,draw,fill=blue!20] (FL) at (-3,8) {};
\node[rounded corners,rectangle,draw,fill=blue!20] (GLR) at (0,8) {};
\node[rounded corners,rectangle,draw,fill=blue!20] (FR) at (3,8) {};
\node[rounded corners,rectangle,draw,fill=blue!20] (PL) at (-3,6) {};
\node[rounded corners,rectangle,draw,fill=blue!20] (F) at (0,6) {};
\node[rounded corners,rectangle,draw,fill=blue!20] (PR) at (3,6) {};
\node[rounded corners,rectangle,draw,fill=blue!20] (GL) at (-3,4) {};
\node[rounded corners,rectangle,draw,fill=blue!20] (P) at (0,4) {};
\node[rounded corners,rectangle,draw,fill=blue!20] (GR) at (3,4) {};
\node[rounded corners,rectangle,draw,fill=blue!20] (G) at (0,2) {};
\node[rounded corners,rectangle,draw,fill=blue!20] (E) at (6,4) {};
\node[rounded corners,rectangle,draw,fill=blue!20] (Q) at (6,2){};
\node[rounded corners,rectangle,draw,fill=blue!20] (1) at (6,0) {};
\draw (FLR)--(FL)--(PL)--(GL)--(G)--(GR)--(PR)--(FR)--(FLR)--(GLR)--(PL)--(P)--(PR)--(GLR) (F)--(P)--(G) (E)--(Q)--(1);
\draw[line width=1.5mm,white] (FL)--(F)--(FR);
\draw (FL)--(F)--(FR);
\draw (FLR)--(GLR);
\draw[line width=1.5mm,white] (P)--(Q) (F)--(E);
\draw (F)--(E) (P)--(Q) (G)--(1);
\node[rounded corners,rectangle,draw,fill=blue!20] (FLR) at (0,10) {\small $M=F_{LR}$};
\node[rounded corners,rectangle,draw,fill=blue!20] (FL) at (-3,8) {\small $F_L$};
\node[rounded corners,rectangle,draw,fill=blue!20] (GLR) at (0,8) {\small $G_{LR}$};
\node[rounded corners,rectangle,draw,fill=blue!20] (FR) at (3,8) {\small $F_R$};
\node[rounded corners,rectangle,draw,fill=blue!20] (PL) at (-3,6) {\small $P_L$};
\node[rounded corners,rectangle,draw,fill=blue!20] (F) at (0,6) {\small $F$};
\node[rounded corners,rectangle,draw,fill=blue!20] (PR) at (3,6) {\small $P_R$};
\node[rounded corners,rectangle,draw,fill=blue!20] (GL) at (-3,4) {\small $G_L$};
\node[rounded corners,rectangle,draw,fill=blue!20] (P) at (0,4) {\small $P$};
\node[rounded corners,rectangle,draw,fill=blue!20] (GR) at (3,4) {\small $G_R$};
\node[rounded corners,rectangle,draw,fill=blue!20] (G) at (0,2) {\small $G$};
\node[rounded corners,rectangle,draw,fill=blue!20] (E) at (6,4) {\small $E$};
\node[rounded corners,rectangle,draw,fill=blue!20] (Q) at (6,2) {\small $Q$};
\node[rounded corners,rectangle,draw,fill=blue!20] (1) at (6,0) {\small $1$};
\node () at (0,13.2) {$(0,1,0,0)$};
\end{scope}
\begin{scope}[shift={(34,0)}]
\node[rounded corners,rectangle,draw,fill=blue!20] (FLR) at (0,10) {};
\node[rounded corners,rectangle,draw,fill=blue!20] (FL) at (-3,8) {};
\node[rounded corners,rectangle,draw,fill=blue!20] (GLR) at (0,8) {};
\node[rounded corners,rectangle,draw,fill=blue!20] (FR) at (3,8) {};
\node[rounded corners,rectangle,draw,fill=blue!20] (PL) at (-3,6) {};
\node[rounded corners,rectangle,draw,fill=blue!20] (F) at (0,6) {};
\node[rounded corners,rectangle,draw,fill=blue!20] (PR) at (3,6) {};
\node[rounded corners,rectangle,draw,fill=blue!20] (GL) at (-3,4) {};
\node[rounded corners,rectangle,draw,fill=blue!20] (P) at (0,4) {};
\node[rounded corners,rectangle,draw,fill=blue!20] (GR) at (3,4) {};
\node[rounded corners,rectangle,draw,fill=blue!20] (G) at (0,2) {};
\draw (FLR)--(FL)--(PL)--(GL)--(G)--(GR)--(PR)--(FR)--(FLR)--(GLR)--(PL)--(P)--(PR)--(GLR) (F)--(P)--(G);
\draw[line width=1.5mm,white] (FL)--(F)--(FR);
\draw (FL)--(F)--(FR);
\draw (FLR)--(GLR);
\node[rounded corners,rectangle,draw,fill=blue!20] (FLR) at (0,10) {\small $M=F_{LR}$};
\node[rounded corners,rectangle,draw,fill=blue!20] (FL) at (-3,8) {\small $F_L$};
\node[rounded corners,rectangle,draw,fill=blue!20] (GLR) at (0,8) {\small $G_{LR}$};
\node[rounded corners,rectangle,draw,fill=blue!20] (FR) at (3,8) {\small $F_R$};
\node[rounded corners,rectangle,draw,fill=blue!20] (PL) at (-3,6) {\small $P_L$};
\node[rounded corners,rectangle,draw,fill=blue!20] (F) at (0,6) {\small $F=E$};
\node[rounded corners,rectangle,draw,fill=blue!20] (PR) at (3,6) {\small $P_R$};
\node[rounded corners,rectangle,draw,fill=blue!20] (GL) at (-3,4) {\small $G_L$};
\node[rounded corners,rectangle,draw,fill=blue!20] (P) at (0,4) {\small $P=Q$};
\node[rounded corners,rectangle,draw,fill=blue!20] (GR) at (3,4) {\small $G_R$};
\node[rounded corners,rectangle,draw,fill=blue!20] (G) at (0,2) {\small $G=1$};
\node () at (0,13.2) {$(0,1,0,1)$};
\end{scope}
\end{tikzpicture}
}
\caption{The lattice $\Lp(M)$ when $M$ has type $\T(M)=(0,i,0,j)$.  In each case, the nodes represent distinct submonoids of $M$.  For other types we have $\Lp(M)=\L(M)$; cf.~Figures \ref{f:L1} and \ref{f:L2}.}
\label{f:L4}
\end{center}
\end{figure}

The exact shape of $\Lp(M)$ depends on the values of $i=\Ttwo(M)$ and $j=\Tfour(M)$, and these determine which thin edges (if any) of Figure \ref{f:lattice_red2} to contract (but keep Lemma \ref{l:PQ}\ref{PQ1} in mind).  The possible shapes are shown in Figure \ref{f:L4}.
Figure \ref{f:L5} shows the possible lattices $\Lp(M)$, for an arbitrary monoid $M$, up to lattice isomorphism (cf.~Figure \ref{f:L3}).

\begin{figure}[ht]
\begin{center}
\begin{tikzpicture}[scale=0.3]
\begin{scope}[shift={(-2,2)}]
\node[circle,fill=black, inner sep = 0.0cm] (FLR) at (0,10) {};
\node[circle,fill=black, inner sep = 0.0cm] (FL) at (-3,8) {};
\node[circle,fill=black, inner sep = 0.0cm] (FR) at (3,8) {};
\node[circle,fill=black, inner sep = 0.0cm] (F) at (0,6) {};
\node[circle,fill=black, inner sep = 0.0cm] (PL) at (-3,6) {};
\node[circle,fill=black, inner sep = 0.0cm] (PR) at (3,6) {};
\node[circle,fill=black, inner sep = 0.0cm] (P) at (0,4) {};
\node[circle,fill=black, inner sep = 0.0cm] (GL) at (-3,4) {};
\node[circle,fill=black, inner sep = 0.0cm] (GR) at (3,4) {};
\node[circle,fill=black, inner sep = 0.0cm] (G) at (0,2) {};
\node[circle,fill=black, inner sep = 0.0cm] (GLR) at (0,8) {}; 
\draw (FLR)--(FL)--(PL)--(GL)--(G)--(GR)--(PR)--(FR)--(FLR)--(GLR)--(PL)--(P)--(PR)--(GLR) (F)--(P)--(G) (FL)--(F)--(FR);
\node[circle,fill=black, inner sep = 0.0cm] (E) at (6,4) {}; \node[circle,fill=black, inner sep = 0.0cm] (Q) at (6,2) {}; \node[circle,fill=black, inner sep = 0.0cm] (O) at (6,0) {}; \draw (F)--(E)--(O)--(G) (P)--(Q);
\node[circle,fill=black, inner sep = 0.0cm] (M) at (0,12) {}; \draw(M)--(FLR);
\node[circle,fill=black, inner sep = 0.06cm] (FLR) at (0,10) {};
\node[circle,fill=black, inner sep = 0.06cm] (FL) at (-3,8) {};
\node[circle,fill=black, inner sep = 0.06cm] (FR) at (3,8) {};
\node[circle,fill=black, inner sep = 0.06cm] (F) at (0,6) {};
\node[circle,fill=black, inner sep = 0.06cm] (PL) at (-3,6) {};
\node[circle,fill=black, inner sep = 0.06cm] (PR) at (3,6) {};
\node[circle,fill=black, inner sep = 0.06cm] (P) at (0,4) {};
\node[circle,fill=black, inner sep = 0.06cm] (GL) at (-3,4) {};
\node[circle,fill=black, inner sep = 0.06cm] (GR) at (3,4) {};
\node[circle,fill=black, inner sep = 0.06cm] (G) at (0,2) {};
\node[circle,fill=black, inner sep = 0.06cm] (GLR) at (0,8) {}; 
\node[circle,fill=black, inner sep = 0.06cm] (M) at (0,12) {}; 
\node[circle,fill=black, inner sep = 0.06cm] (E) at (6,4) {}; \node[circle,fill=black, inner sep = 0.06cm] (Q) at (6,2) {}; \node[circle,fill=black, inner sep = 0.06cm] (O) at (6,0) {}; 
\end{scope}
\begin{scope}[shift={(10,2)}]
\node[circle,fill=black, inner sep = 0.0cm] (FLR) at (0,10) {};
\node[circle,fill=black, inner sep = 0.0cm] (FL) at (-3,8) {};
\node[circle,fill=black, inner sep = 0.0cm] (FR) at (3,8) {};
\node[circle,fill=black, inner sep = 0.0cm] (F) at (0,6) {};
\node[circle,fill=black, inner sep = 0.0cm] (PL) at (-3,6) {};
\node[circle,fill=black, inner sep = 0.0cm] (PR) at (3,6) {};
\node[circle,fill=black, inner sep = 0.0cm] (P) at (0,4) {};
\node[circle,fill=black, inner sep = 0.0cm] (GL) at (-3,4) {};
\node[circle,fill=black, inner sep = 0.0cm] (GR) at (3,4) {};
\node[circle,fill=black, inner sep = 0.0cm] (G) at (0,2) {};
\node[circle,fill=black, inner sep = 0.0cm] (GLR) at (0,8) {}; 
\draw (FLR)--(FL)--(PL)--(GL)--(G)--(GR)--(PR)--(FR)--(FLR)--(GLR)--(PL)--(P)--(PR)--(GLR) (F)--(P)--(G) (FL)--(F)--(FR);
\node[circle,fill=black, inner sep = 0.0cm] (M) at (0,12) {}; \draw(M)--(FLR);
\node[circle,fill=black, inner sep = 0.06cm] (FLR) at (0,10) {};
\node[circle,fill=black, inner sep = 0.06cm] (FL) at (-3,8) {};
\node[circle,fill=black, inner sep = 0.06cm] (FR) at (3,8) {};
\node[circle,fill=black, inner sep = 0.06cm] (F) at (0,6) {};
\node[circle,fill=black, inner sep = 0.06cm] (PL) at (-3,6) {};
\node[circle,fill=black, inner sep = 0.06cm] (PR) at (3,6) {};
\node[circle,fill=black, inner sep = 0.06cm] (P) at (0,4) {};
\node[circle,fill=black, inner sep = 0.06cm] (GL) at (-3,4) {};
\node[circle,fill=black, inner sep = 0.06cm] (GR) at (3,4) {};
\node[circle,fill=black, inner sep = 0.06cm] (G) at (0,2) {};
\node[circle,fill=black, inner sep = 0.06cm] (GLR) at (0,8) {}; 
\node[circle,fill=black, inner sep = 0.06cm] (M) at (0,12) {}; 
\end{scope}
\begin{scope}[shift={(20,2)}]
\node[circle,fill=black, inner sep = 0.0cm] (FLR) at (0,8) {};
\node[circle,fill=black, inner sep = 0.0cm] (FL) at (-3,6) {};
\node[circle,fill=black, inner sep = 0.0cm] (FR) at (3,6) {};
\node[circle,fill=black, inner sep = 0.0cm] (F) at (0,4) {};
\node[circle,fill=black, inner sep = 0.0cm] (GL) at (-3,4) {};
\node[circle,fill=black, inner sep = 0.0cm] (GR) at (3,4) {};
\node[circle,fill=black, inner sep = 0.0cm] (G) at (0,2) {};
\draw (G)--(GR)--(FR)--(FLR)--(FL)--(GL)--(G) (FL)--(F)--(FR) (G)--(F);
\node[circle,fill=black, inner sep = 0.0cm] (E) at (3,2) {}; \node[circle,fill=black, inner sep = 0.0cm] (O) at (3,0) {}; \draw (F)--(E)--(O)--(G);
\node[circle,fill=black, inner sep = 0.0cm] (M) at (0,10) {}; \draw(M)--(FLR);
\node[circle,fill=black, inner sep = 0.06cm] (FLR) at (0,8) {};
\node[circle,fill=black, inner sep = 0.06cm] (FL) at (-3,6) {};
\node[circle,fill=black, inner sep = 0.06cm] (FR) at (3,6) {};
\node[circle,fill=black, inner sep = 0.06cm] (F) at (0,4) {};
\node[circle,fill=black, inner sep = 0.06cm] (GL) at (-3,4) {};
\node[circle,fill=black, inner sep = 0.06cm] (GR) at (3,4) {};
\node[circle,fill=black, inner sep = 0.06cm] (G) at (0,2) {};
\node[circle,fill=black, inner sep = 0.06cm] (E) at (3,2) {}; \node[circle,fill=black, inner sep = 0.06cm] (O) at (3,0) {}; 
\node[circle,fill=black, inner sep = 0.06cm] (M) at (0,10) {}; 
\end{scope}
\begin{scope}[shift={(30,2)}]
\node[circle,fill=black, inner sep = 0.0cm] (FLR) at (0,8) {};
\node[circle,fill=black, inner sep = 0.0cm] (FL) at (-3,6) {};
\node[circle,fill=black, inner sep = 0.0cm] (FR) at (3,6) {};
\node[circle,fill=black, inner sep = 0.0cm] (F) at (0,4) {};
\node[circle,fill=black, inner sep = 0.0cm] (GL) at (-3,4) {};
\node[circle,fill=black, inner sep = 0.0cm] (GR) at (3,4) {};
\node[circle,fill=black, inner sep = 0.0cm] (G) at (0,2) {};
\draw (G)--(GR)--(FR)--(FLR)--(FL)--(GL)--(G) (FL)--(F)--(FR) (G)--(F);
\node[circle,fill=black, inner sep = 0.0cm] (M) at (0,10) {}; \draw(M)--(FLR);
\node[circle,fill=black, inner sep = 0.06cm] (FLR) at (0,8) {};
\node[circle,fill=black, inner sep = 0.06cm] (FL) at (-3,6) {};
\node[circle,fill=black, inner sep = 0.06cm] (FR) at (3,6) {};
\node[circle,fill=black, inner sep = 0.06cm] (F) at (0,4) {};
\node[circle,fill=black, inner sep = 0.06cm] (GL) at (-3,4) {};
\node[circle,fill=black, inner sep = 0.06cm] (GR) at (3,4) {};
\node[circle,fill=black, inner sep = 0.06cm] (G) at (0,2) {};
\node[circle,fill=black, inner sep = 0.06cm] (M) at (0,10) {}; 
\end{scope}
\begin{scope}[shift={(-2,-10)}]
\node[circle,fill=black, inner sep = 0.0cm] (FLR) at (0,10) {};
\node[circle,fill=black, inner sep = 0.0cm] (FL) at (-3,8) {};
\node[circle,fill=black, inner sep = 0.0cm] (FR) at (3,8) {};
\node[circle,fill=black, inner sep = 0.0cm] (F) at (0,6) {};
\node[circle,fill=black, inner sep = 0.0cm] (PL) at (-3,6) {};
\node[circle,fill=black, inner sep = 0.0cm] (PR) at (3,6) {};
\node[circle,fill=black, inner sep = 0.0cm] (P) at (0,4) {};
\node[circle,fill=black, inner sep = 0.0cm] (GL) at (-3,4) {};
\node[circle,fill=black, inner sep = 0.0cm] (GR) at (3,4) {};
\node[circle,fill=black, inner sep = 0.0cm] (G) at (0,2) {};
\node[circle,fill=black, inner sep = 0.0cm] (GLR) at (0,8) {}; 
\draw (FLR)--(FL)--(PL)--(GL)--(G)--(GR)--(PR)--(FR)--(FLR)--(GLR)--(PL)--(P)--(PR)--(GLR) (F)--(P)--(G) (FL)--(F)--(FR);
\node[circle,fill=black, inner sep = 0.0cm] (E) at (6,4) {}; \node[circle,fill=black, inner sep = 0.0cm] (Q) at (6,2) {}; \node[circle,fill=black, inner sep = 0.0cm] (O) at (6,0) {}; \draw (F)--(E)--(O)--(G) (P)--(Q);
\node[circle,fill=black, inner sep = 0.06cm] (FLR) at (0,10) {};
\node[circle,fill=black, inner sep = 0.06cm] (FL) at (-3,8) {};
\node[circle,fill=black, inner sep = 0.06cm] (FR) at (3,8) {};
\node[circle,fill=black, inner sep = 0.06cm] (F) at (0,6) {};
\node[circle,fill=black, inner sep = 0.06cm] (PL) at (-3,6) {};
\node[circle,fill=black, inner sep = 0.06cm] (PR) at (3,6) {};
\node[circle,fill=black, inner sep = 0.06cm] (P) at (0,4) {};
\node[circle,fill=black, inner sep = 0.06cm] (GL) at (-3,4) {};
\node[circle,fill=black, inner sep = 0.06cm] (GR) at (3,4) {};
\node[circle,fill=black, inner sep = 0.06cm] (G) at (0,2) {};
\node[circle,fill=black, inner sep = 0.06cm] (GLR) at (0,8) {}; 
\node[circle,fill=black, inner sep = 0.06cm] (E) at (6,4) {}; \node[circle,fill=black, inner sep = 0.06cm] (Q) at (6,2) {}; \node[circle,fill=black, inner sep = 0.06cm] (O) at (6,0) {}; 
\end{scope}
\begin{scope}[shift={(10,-10)}]
\node[circle,fill=black, inner sep = 0.0cm] (FLR) at (0,10) {};
\node[circle,fill=black, inner sep = 0.0cm] (FL) at (-3,8) {};
\node[circle,fill=black, inner sep = 0.0cm] (FR) at (3,8) {};
\node[circle,fill=black, inner sep = 0.0cm] (F) at (0,6) {};
\node[circle,fill=black, inner sep = 0.0cm] (PL) at (-3,6) {};
\node[circle,fill=black, inner sep = 0.0cm] (PR) at (3,6) {};
\node[circle,fill=black, inner sep = 0.0cm] (P) at (0,4) {};
\node[circle,fill=black, inner sep = 0.0cm] (GL) at (-3,4) {};
\node[circle,fill=black, inner sep = 0.0cm] (GR) at (3,4) {};
\node[circle,fill=black, inner sep = 0.0cm] (G) at (0,2) {};
\node[circle,fill=black, inner sep = 0.0cm] (GLR) at (0,8) {}; 
\draw (FLR)--(FL)--(PL)--(GL)--(G)--(GR)--(PR)--(FR)--(FLR)--(GLR)--(PL)--(P)--(PR)--(GLR) (F)--(P)--(G) (FL)--(F)--(FR);
\node[circle,fill=black, inner sep = 0.06cm] (FLR) at (0,10) {};
\node[circle,fill=black, inner sep = 0.06cm] (FL) at (-3,8) {};
\node[circle,fill=black, inner sep = 0.06cm] (FR) at (3,8) {};
\node[circle,fill=black, inner sep = 0.06cm] (F) at (0,6) {};
\node[circle,fill=black, inner sep = 0.06cm] (PL) at (-3,6) {};
\node[circle,fill=black, inner sep = 0.06cm] (PR) at (3,6) {};
\node[circle,fill=black, inner sep = 0.06cm] (P) at (0,4) {};
\node[circle,fill=black, inner sep = 0.06cm] (GL) at (-3,4) {};
\node[circle,fill=black, inner sep = 0.06cm] (GR) at (3,4) {};
\node[circle,fill=black, inner sep = 0.06cm] (G) at (0,2) {};
\node[circle,fill=black, inner sep = 0.06cm] (GLR) at (0,8) {}; 
\end{scope}
\begin{scope}[shift={(20,-10)}]
\node[circle,fill=black, inner sep = 0.0cm] (FLR) at (0,8) {};
\node[circle,fill=black, inner sep = 0.0cm] (FL) at (-3,6) {};
\node[circle,fill=black, inner sep = 0.0cm] (FR) at (3,6) {};
\node[circle,fill=black, inner sep = 0.0cm] (F) at (0,4) {};
\node[circle,fill=black, inner sep = 0.0cm] (GL) at (-3,4) {};
\node[circle,fill=black, inner sep = 0.0cm] (GR) at (3,4) {};
\node[circle,fill=black, inner sep = 0.0cm] (G) at (0,2) {};
\draw (G)--(GR)--(FR)--(FLR)--(FL)--(GL)--(G) (FL)--(F)--(FR) (G)--(F);
\node[circle,fill=black, inner sep = 0.0cm] (E) at (3,2) {}; \node[circle,fill=black, inner sep = 0.0cm] (O) at (3,0) {}; \draw (F)--(E)--(O)--(G);
\node[circle,fill=black, inner sep = 0.06cm] (FLR) at (0,8) {};
\node[circle,fill=black, inner sep = 0.06cm] (FL) at (-3,6) {};
\node[circle,fill=black, inner sep = 0.06cm] (FR) at (3,6) {};
\node[circle,fill=black, inner sep = 0.06cm] (F) at (0,4) {};
\node[circle,fill=black, inner sep = 0.06cm] (GL) at (-3,4) {};
\node[circle,fill=black, inner sep = 0.06cm] (GR) at (3,4) {};
\node[circle,fill=black, inner sep = 0.06cm] (G) at (0,2) {};
\node[circle,fill=black, inner sep = 0.06cm] (E) at (3,2) {}; \node[circle,fill=black, inner sep = 0.06cm] (O) at (3,0) {}; 
\end{scope}
\begin{scope}[shift={(30,-10)}]
\node[circle,fill=black, inner sep = 0.0cm] (FLR) at (0,8) {};
\node[circle,fill=black, inner sep = 0.0cm] (FL) at (-3,6) {};
\node[circle,fill=black, inner sep = 0.0cm] (FR) at (3,6) {};
\node[circle,fill=black, inner sep = 0.0cm] (F) at (0,4) {};
\node[circle,fill=black, inner sep = 0.0cm] (GL) at (-3,4) {};
\node[circle,fill=black, inner sep = 0.0cm] (GR) at (3,4) {};
\node[circle,fill=black, inner sep = 0.0cm] (G) at (0,2) {};
\draw (G)--(GR)--(FR)--(FLR)--(FL)--(GL)--(G) (FL)--(F)--(FR) (G)--(F);
\node[circle,fill=black, inner sep = 0.06cm] (FLR) at (0,8) {};
\node[circle,fill=black, inner sep = 0.06cm] (FL) at (-3,6) {};
\node[circle,fill=black, inner sep = 0.06cm] (FR) at (3,6) {};
\node[circle,fill=black, inner sep = 0.06cm] (F) at (0,4) {};
\node[circle,fill=black, inner sep = 0.06cm] (GL) at (-3,4) {};
\node[circle,fill=black, inner sep = 0.06cm] (GR) at (3,4) {};
\node[circle,fill=black, inner sep = 0.06cm] (G) at (0,2) {};
\end{scope}
\begin{scope}[shift={(-2,-22)}]
\node[circle,fill=black, inner sep = 0.0cm] (M) at (0,10) {};
\node[circle,fill=black, inner sep = 0.0cm] (F) at (0,8) {};
\node[circle,fill=black, inner sep = 0.0cm] (G) at (-2,6) {};
\node[circle,fill=black, inner sep = 0.0cm] (E) at (2,6) {};
\node[circle,fill=black, inner sep = 0.0cm] (1) at (0,4) {};
\draw (M)--(F);
\draw (F)--(E) (G)--(1);
\draw (F)--(G) (E)--(1);
\node[circle,fill=black, inner sep = 0.06cm] (M) at (0,10) {};
\node[circle,fill=black, inner sep = 0.06cm] (F) at (0,8) {};
\node[circle,fill=black, inner sep = 0.06cm] (G) at (-2,6) {};
\node[circle,fill=black, inner sep = 0.06cm] (E) at (2,6) {};
\node[circle,fill=black, inner sep = 0.06cm] (1) at (0,4) {};
\end{scope}
\begin{scope}[shift={(10,-22)}]
\node[circle,fill=black, inner sep = 0.0cm] (F) at (0,8) {};
\node[circle,fill=black, inner sep = 0.0cm] (G) at (-2,6) {};
\node[circle,fill=black, inner sep = 0.0cm] (E) at (2,6) {};
\node[circle,fill=black, inner sep = 0.0cm] (1) at (0,4) {};
\draw (F)--(E) (G)--(1);
\draw (F)--(G) (E)--(1);
\node[circle,fill=black, inner sep = 0.06cm] (F) at (0,8) {};
\node[circle,fill=black, inner sep = 0.06cm] (G) at (-2,6) {};
\node[circle,fill=black, inner sep = 0.06cm] (E) at (2,6) {};
\node[circle,fill=black, inner sep = 0.06cm] (1) at (0,4) {};
\end{scope}
\begin{scope}[shift={(20,-22)}]
\node[circle,fill=black, inner sep = 0.0cm] (M) at (0,8) {};
\node[circle,fill=black, inner sep = 0.0cm] (F) at (0,6) {};
\node[circle,fill=black, inner sep = 0.0cm] (1) at (0,4) {};
\draw (M)--(F);
\draw (F)--(1);
\node[circle,fill=black, inner sep = 0.06cm] (M) at (0,8) {};
\node[circle,fill=black, inner sep = 0.06cm] (F) at (0,6) {};
\node[circle,fill=black, inner sep = 0.06cm] (1) at (0,4) {};
\end{scope}
\begin{scope}[shift={(27,-22)}]
\node[circle,fill=black, inner sep = 0.0cm] (F) at (0,6) {};
\node[circle,fill=black, inner sep = 0.0cm] (1) at (0,4) {};
\draw (F)--(1);
\node[circle,fill=black, inner sep = 0.06cm] (F) at (0,6) {};
\node[circle,fill=black, inner sep = 0.06cm] (1) at (0,4) {};
\end{scope}
\begin{scope}[shift={(33,-22)}]
\node[circle,fill=black, inner sep = 0.06cm] (1) at (0,4) {};
\end{scope}
\end{tikzpicture}
\caption{The possible lattices $\Lp(M)$ for a monoid $M$, up to lattice isomorphism.}
\label{f:L5}
\end{center}
\end{figure}

\begin{rem}
Recall that $B^0$ is the bicyclic monoid with a zero adjoined.  We noted in Remark \ref{r:notlattice2} that~$\L(B^0)$ is not a sublattice of $\Sub(B^0)$, citing the fact that $F_L\cap G_{LR}\not\in\L(B^0)$, using the usual abbreviations.  However, consulting Tables \ref{t:GNB_small} and \ref{t:GNB_smaller}, we see that $F_L\cap G_{LR}=P_L\in\Lp(B^0)$.  In fact, $\Lp(B^0)$ is closed under arbitrary intersections, as one may easily check using the aforementioned tables, which means that $\Lp(B^0)$ is a sublattice of $\Sub(B^0)$.  The author does not currently know if $\Lp(M)$ is a sublattice of~$\Sub(M)$ for an arbitrary monoid $M$.
\end{rem}

\footnotesize
\def\bibspacing{-1.1pt}
\bibliography{biblio}
\bibliographystyle{abbrv}
\end{document}